\documentclass[journal]{IEEEtran}
\usepackage{nicefrac}
\usepackage{verbatim}
\usepackage{microtype}
\usepackage{amsmath,amssymb,amsfonts,amsthm,amscd,bm}
\usepackage[ruled,linesnumbered]{algorithm2e}
\usepackage{cite}
\usepackage{graphicx}
\usepackage{multirow}
\usepackage{colortbl}
\usepackage{array,booktabs}
\usepackage{makecell}
\usepackage{pifont}

\usepackage{etoolbox}
\makeatletter
\patchcmd{\@makecaption}
  {\scshape}
  {}
  {}
  {}
\makeatother

\newtheorem{theorem}{Theorem}

\newtheorem{assumption}[theorem]{Assumption}
\newtheorem{lemma}[theorem]{Lemma}

\DeclareMathOperator*{\argmin}{arg\,min}

\DeclareMathOperator{\diag}{diag}

\newcommand{\G}{\mathcal{G}}

\newcommand{\E}{\mathcal{E}}
\newcommand{\V}{\mathcal{V}}

\newcommand{\hV}{\hat{\mathcal{V}}}
\newcommand{\D}{\mathcal{D}}
\newcommand{\W}{\mathcal{W}}

\renewcommand{\AA}{\mathcal{A}}

\newcommand{\LL}{\mathcal{L}}

\newcommand{\hM}{\hat{M}}
\newcommand{\MC}{\mathcal{M}}
\newcommand{\hMC}{\hat{\mathcal{M}}}
\newcommand{\HH}{\mathcal{H}}
\newcommand{\hHH}{\hat{\mathcal{H}}}
\newcommand{\vpi}{\bm{\pi}}
\newcommand{\hvpi}{\hat{\bm{\pi}}}

\newcommand{\vx}{\bm{x}}
\newcommand{\vz}{\bm{z}}

\newcommand{\vu}{\bm{u}}

\newcommand{\vn}{\bm{n}}

\newcommand{\vv}{\bm{v}}

\newcommand{\bomega}{\bar{\omega}}

\begin{document}
\title{Distributed Optimization,  Averaging via ADMM, \\ and Network Topology}
\author{Guilherme Fran\c ca, and Jos\' e Bento%
\thanks{GF is with the Mathematical Institute for Data Science,
Johns Hopkins University, Baltimore, MD.
e-mail: guifranca@gmail.com}
\thanks{JB is with the Department of Computer Science,
Boston College, Boston, MA.
e-mail: jose.bento@bc.edu}
}

\maketitle

\begin{abstract}
There has been an increasing necessity for scalable optimization methods, especially due to the explosion in the size of datasets and model complexity in modern machine learning applications.
Scalable solvers often distribute the computation
over a network of processing units.
For simple algorithms such as gradient descent the dependency of the convergence time with the topology of this network is well-known.
However, for more involved algorithms such as the Alternating Direction Methods of Multipliers (ADMM) much less is known.
At the heart of many distributed optimization algorithms there exists a
gossip subroutine which averages local information over the network, and whose efficiency is crucial for the overall
performance of the method.
In this paper
we review recent research in this area and, with the goal
of isolating such a communication exchange behaviour, we
compare different algorithms when applied to a
canonical distributed averaging consensus problem.
We also show interesting connections between ADMM and lifted Markov chains besides providing an explicitly characterization of its convergence and optimal parameter tuning in terms of spectral properties of the network.
Finally, we empirically study the connection between network topology and convergence rates for different algorithms on a real world problem of sensor localization.
\end{abstract}


\IEEEpeerreviewmaketitle

\section{Introduction} \label{sec:introduction}

\IEEEPARstart{D}{uring} the last decade there has been a fast growth on the
size of datasets and complexity of the models used in machine learning and statistics.
At the heart of many applications an optimization problem
needs to be solved, driving the need for efficient solvers whose computations can be parallelized and distributed.
Gradient Descent (GD)---in its many distributed and stochastic variants---is a popular algorithm specially because
its simplicity helps scalable implementations.

When differentiability assumptions of gradient based methods are restrictive a powerful alternative are methods that 
use \emph{proximal maps}, which are often well-defined despite nonsmoothness and even nonconvexity
of the objective function and constraints. One example of a proximal based algorithm is the Alternating Direction
Method of Multipliers (ADMM).
ADMM was
proposed in the 70's \cite{Gabay:1976,Glowinsky:1975,RockaF} but in the last decade has regained great attention \cite{Boyd:2010,Eckstein} thanks to its robustness, extremely mild assumptions for convergence,
and because it can be easily distributed.
Empirically, 
ADMM works well for many hard optimization problems 
while many other methods struggle to converge or do not scale \cite{Boyd:2010,Candes:2011,Derbinski1,Bento1,Bento2,Krishnan,Bento3,yang2019estimating,arminstratis2020}. Recently, there has been promising theoretical results for nonconvex
problems as well \cite{Hong:2016,Yin1}.
ADMM is known to have fast converge when a low-accuracy solution is the goal, but on the other hand its
convergence rate can be sensitive to parameter tuning.
Although empirically ADMM has great performance, 
on a theoretical level its convergence and stability properties are
still not fully understood, specially in distributed settings. For example, only recently \cite{FrancaBento,giselsson2017linear,Jordan} explicit---i.e., with no unspecified constants---and tight bounds on its convergence rate have been proven. 
Interesting connections of ADMM, and accelerated variants thereof,  with  continuous dynamical
systems have also been recently proposed \cite{Franca:2018_admm,Franca:2018_nonsmooth,Franca:2019_splitting}.

An important class of problems in distributed optimization is the class of
\emph{consensus problems}:
\begin{equation}
\label{eq:consensus}
\min_{\vz \in \mathbb{R}^p} \bigg\{ f(\vz) = \sum_{i=1}^{n} f_i(\vz) \bigg\},
\end{equation}
where $f_i : \mathbb{R}^p \to \mathbb{R}$ is one term of the
objective function $f$.
A distributed  solver typically
consists of $n$ ``agents,'' each being a computation unit.
Agent $i$ only accesses
$f_i$ and has a local copy of (some few components of) $\vz$, which can be high-dimensional. Each agent communicates only with a subset of neighboring agents according to a graph
$\G = (\V, \E)$, where  $\V \equiv \{1,\dots,n\}$ is the set of nodes and $\E$ is the set of edges. Node $i\in \V$ represents agent $i$ and edge $(i,j)\in \E$ is a communication link between agents $i$ and $j$.
The goal is that the ``opinion'' of different agents, i.e., their local copies of $\vz$, converge to an agreement so that a solution
to \eqref{eq:consensus} is obtained.

To design methods
that respect this communication and locality constraints, one often rewrites \eqref{eq:consensus} in the form
\begin{equation}
\label{eq:consensus_distributed}
\min_{\{\vx_i \in \mathbb{R}^p\}_{i \in \V}}  \sum_{i \in \V} f_i(\vx_i) + \sum_{(i,j) \in \E} \mathbb{I}(\vx_i,\vx_j),
\end{equation}
where $\mathbb{I}(\vx_i,\vx_j)$ takes the value of $0$ if $\vx_i  = \vx_j$ and $\infty$ otherwise.
In this form, during each iteration of a distributed method, agent $i$ updates its local copy, $\vx_i$, to decrease the value of its local function $f_i$, while simultaneously communicating with its neighbours to enforce local constraints $\vx_i = \vx_j$ for all $j$'s in its neighborhood.

For a nice review on the convergence of distributed algorithms with the topology of $\G$, see \cite{nedic2018network}. The first part of this paper considers distributed averaging via $\vz^{t+1}  = A^t \vz^t$, where $t$ is the iteration time, $\vz^t \equiv \{z^t_i\}_{i \in \V}$ is distributed across $\V$, and the gossip matrix $A^t$---assumed to be doubly stochastic---models message exchange along $\E$. In the second part, these results for averaging are employed to study algorithms for minimizing $f(\vz)$ via $\vz^{t+1} = P(A^{t} \vz^t - \bm{g}^t)$, where $\bm{g}^t$ is some gradient information about $f$ and $P$ is a projection operator.
We note that \cite{nedic2018network} also considers time-varying graphs, which will not be considered in this paper.
Unfortunately, the algorithms considered in \cite{nedic2018network} do not link convergence time with topology in an optimal manner. More specifically, let $\omega_{n}$ be the normalized eigengap of the Lapacian, $L$, of $\G$, i.e., the ratio between the second smallest and the largest eigenvalues:
\begin{equation}\label{eq:spectral_gap}
{\omega}_{n} \equiv (\lambda_{|\V| -1 }(L))/(\lambda_{1}(L)).
\end{equation}
The convergence time of the algorithms in \cite{nedic2018network}  is of $\mathcal{O}(1/{\omega}_{n})$.
However, there are known algorithms with convergence time of $\mathcal{O}(1/\sqrt{{\omega}_{n}})$. Part of this suboptimality is due to the requirement that $A^t$ is (doubly) stochastic.
Some works \cite{makhdoumi2014broadcast,makhdoumi2017convergence}  relate convergence time and topology when solving \eqref{eq:consensus} for ADMM-type algorithms; e.g., \cite{makhdoumi2017convergence} obtains a convergence time of
 $\mathcal{O}(1/ \omega_n )$, which is again suboptimal.
Classical results on distributed optimization  \cite{BertsekasBook} are either too general or suboptimal with respect to the topology of $\G$.

Recently, algorithms whose convergence time depends \emph{optimally}
on the topology of $\G$, i.e. $\mathcal{O}(1/\sqrt{\omega_n})$,  have been proposed  \cite{seaman2017optimal,scaman2018optimal,XFilter};
moreover, they put forward problems that cannot be solved to fixed accuracy in less than $\Omega(1/\sqrt{\omega_n})$ iterations.
Similarly to \cite{nedic2018network}, such algorithms solve \eqref{eq:consensus} by computing  fixed-point iterations of the first-order optimality conditions by averaging gradient information over $\G$. However, these methods achieve optimality by using an improved distributed averaging or ``gossip'' step;
in particular, \cite{seaman2017optimal,scaman2018optimal,XFilter} use Chebyshev filtering to manipulate the spectrum of the gossip matrix, an idea first used in the context of distributed averaging in \cite{kokiopoulou2008polynomial}.
Another idea that was also explored to improve averaging is lifted Markov chains \cite{4557655,jung2009distributed}, although these works  do not approach the averaging problem from an optimization perspective.
An usual take on distributed averaging, the aforementioned 
bottleneck of many  solvers, is via the optimization problem
\begin{equation}
\min_{\vz \in \mathbb{R}^n}  \sum_{i \in \V} (z_i - c_i)^2 \ \text{ subject to } z_i = z_j \text{ if } (i,j) \in \E, \label{eq:alternative_formulation_for_avg_concensus}
\end{equation}
where $c_i$ are given values intended to be averaged.
The analysis of (a slightly different formulation of) ADMM for this problem has been done  \cite{ghadimi2015optimal,Ghadimi2014averaging}, where optimal tuning rules are provided.
In this current paper, we complement such results by reviewing a nonstandard ADMM formulation for distributed averaging. More specifically, instead of \eqref{eq:alternative_formulation_for_avg_concensus} we consider
\begin{equation} \label{eq:gen_quad}
\min_{ \vz \in \mathbb{R}^n} \bigg \{
f(\bm{z}) = \tfrac{1}{2} \sum_{(i,j)  \in \E} (z_i - z_j)^2
\bigg \},
\end{equation}
which we will refer to as the \emph{canonical consensus problem}.

Deriving distributed averaging algorithms from \eqref{eq:gen_quad}
is  different from the much more common approach through \eqref{eq:alternative_formulation_for_avg_concensus}.
Our reason for doing so is twofold. First, it can be shown that within this formulation, ADMM can be seen as a \emph{lifted} Markov chain \cite{FrancaBentoMarkov}, which achieves such an optimal ``square root speedup'' over GD-type approaches. Similarly to
\cite{nedic2018network}, the resulting algorithm has the form $\vz^{t+1} = T_A \vz^{t}$, however the gossip matrix $T_A$ is not  (doubly) stochastic, allowing one to escape suboptimal bounds \cite{nedic2018network}. Second, this formulation is amenable to an analytical treatment regarding the relationship between convergence time and graph topology  in an explicit manner \cite{FrancaBentoTopology,FrancaBentoTopology2}, yielding usually faster rates than previously obtained for the standard ADMM formulation \cite{Ghadimi2014averaging,ghadimi2015optimal}.

This paper focuses on  comparing how different algorithms 
tackle the distributed averaging step, and in particular their performance in solving  \eqref{eq:gen_quad}. More specifically, in Section~\ref{sec:notation} we collect
basic definitions and introduce notation.
In Section~\ref{sec:distributed_algos}, we review recently proposed distributed optimization algorithms.
In Section~\ref{sec:canonical}, we review interesting connections between ADMM, GD, spectral graph theory, and lifted Markov chains for solving the crucial step of distributed averaging which is common to most optimal (but also non-optimal) general solvers.
We  present explicit formulas for the convergence rate of ADMM and optimal parameters in terms of the spectral properties  of $\G$.
%
%
In Section~5, we show numerical experiments comparing the performance of the algorithms discussed in Section~\ref{sec:distributed_algos}, contrasting with theoretical results. 
We remark that this paper  comes accompanied with publicly available code \cite{GitLink}, where all the algorithms here considered are implemented.\footnote{This can be a useful starting point for practitioners interested in using these distributed optimization methods.}

\section{Preliminaries and Notation}
\label{sec:notation}

In this section we introduce some necessary elements of graph theory; we refer to \cite{GraphTheoryBook1,GraphTheoryBook2} for  details.

Let $\G = (\V, \E)$ be an undirected graph,
with vertex set $\V$ of size $n = |\V|$ and edge set $\E$ of size $m = |\E|$.
We make the following general assumption throughout.
\begin{assumption}\label{thm:assumption}
$\G$ is undirected, connected, and simple.
\end{assumption}
We denote its Laplacian by $L = \D - \AA$, where $\D$ is the degree matrix
and $\AA$ is the adjacency
matrix. 
%
The normalized Laplacian is
$\LL \equiv \D^{-1/2} L \D^{-1/2} = I - \D^{-1/2} \AA \D^{-1/2}$,
and
\begin{equation}\label{eq:def_walk_matrix}
\W \equiv \D^{-1} \AA
\end{equation}
is the probability transition matrix of a random
walk on $\G$ where each node has equal probability to jump to any neighboring node.
Although general stochastic matrices may have complex eigenvalues, in our case
$\W$ is co-spectral to $D^{1/2}\W D^{-1/2} = \D^{-1/2} \AA \D^{-1/2}$, thus, as $L$ and $\LL$, it only has real eigenvalues and is diagonalizable. Its eigenvalues are in the range $[-1,1]$.
%
We denote neighbours of $i \in \V$ by
\begin{equation}
N_i \equiv \{ j \in \V \, | \, (i,j) \in \E  \},
\end{equation}
and the neighbors of $e = (i,j)\in \E$ by $N_e = \{i,j\}$. When clear from the context, $N_i$ might also mean $N_i \equiv \{e\in \E| e = (i,j)\}$, i.e., the set of edges adjacent to node $i$.
The degree of node $i \in \V$ is denoted by $d_i$, and the maximum and minimum degrees by $d_{\max}$ and $d_{\min}$, respectively.

For any matrix $T \in \mathbb{R}^{n\times n}$, unless stated otherwise, we assume that its eigenvalues are ordered by
magnitude:
\begin{equation}\label{eq:ordering}
  |\lambda_1(T)| \ge |\lambda_2(T)| \ge \dotsm \ge |\lambda_n(T)|.
\end{equation}
In particular, if $T$ has only nonnegative eigenvalues, e.g., it is positive semidefinite such as the Laplacian $L$,
we  drop the absolute values above. We refer to $\lambda(T)$ as a generic eigenvalue
when the above ordering is not important in the discussion.

The following problem,
of which \eqref{eq:gen_quad} is a special instance,
will  be used in Section \ref{sec:distributed_algos}:
\begin{equation}\label{eq:pair_interaction_problem}
\min_{ \vz \in \mathbb{R}^n} \bigg \{
f(\bm{z}) = \sum_{e = (i,j) \in \E } f_e(z_i, z_j) \bigg \}.
\end{equation}

A matrix that will play an important role 
is the row stochastic matrix
$S \in \mathbb{R}^{2|\E|\times |\V|}$ defined as follows. Define the extended set of edges:
\begin{equation} \label{eq:extended_E}
\hat{\E} \equiv \{ (e, i): e = (i,j) \in \E \}.
\end{equation}
Denoting $e_i = (e, i)$, we thus have:
\begin{equation}\label{eq:Sdef}
S_{e_i, i} = S_{e_j,j} \equiv \begin{cases}
1 & \mbox{if $e=(i,j) \in \E$,} \\
0 &  \mbox{otherwise.}
\end{cases}
\end{equation}
Note that $S^T S = \D$ (where $S^T$ denotes the transpose of $S$) is the degree matrix of the  graph $\G$.

In general, given an index set $K$ and a vector $\vx$, we let $\vx_{\bm{K}}$ denote a vector containing only those components of $\vx$ with index in $K$. For instance, if $e = (i,j) \in {\E}$ then
\begin{equation} \label{eq:zNe}
\vz_{\bm{N_e}} \equiv (z_i, z_j)^T.
\end{equation}
%

To describe the behaviour of several algorithms with respect to the topology of the underlying graph, it is useful to introduce the following terminology:
\begin{itemize}
    \item $\omega^\star$ is the 2nd largest (not in magnitude) eigenvalue of $\W$;
    \item $\bar{\omega}$ is the 2nd smallest (not in magnitude) eigenvalue of $\W$ that is different than $-1$;
    \item $\hat{\omega}$ is the largest eigenvalue (in magnitude) of $\W$ that different than $1$;
    \item $\hat{\omega}_{\delta} \equiv 1 - \hat{\omega} > 0$;
    \item $\omega_{n}$ is the ratio between smallest nonzero and largest eigenvalues of $L$;
    \item $\omega_L$  is the second smallest eigenvalue of the Laplacian $L$.
\end{itemize}
Such quantities can be computed offline, using standard numerical linear algebra, but also with distributed algorithms; see e.g. \cite{tran2014distributed,di2014distributed,gusrialdi2017distributed}.
For many graphs we have that $\omega^* = \hat{\omega}$, and if  $d_{\min},d_{\max}$ are  bounded, then  $\hat{\omega}_{\delta} = \Theta(\omega_n) = \Theta(\omega_L)$. The latter point follows from the following inequalities \cite{zumstein2005comparison}:
\begin{equation} \label{eq:relation_among_spectral_measures}
\begin{split}
\lambda(\W) &= 1 - \lambda(\mathcal{L}) ,  \\
d_{\textnormal{min}} \lambda_i (\mathcal{L}) &\le \lambda_i(L) \le
d_{\textnormal{max}} \lambda_i (\mathcal{L}) \quad (i=1,\dotsc,|\V|), \\
d_{\textnormal{max}} & \le \lambda_1(L) \le 2 d_{\textnormal{max}}.
\end{split}
\end{equation}%
In particular, usually
$1-\omega^*, \hat{\omega}_{\delta} ,\omega_n$, and $\omega_L$ all converge to zero at the same rate.


We define the \emph{asymptotic convergence
rate}, $\tau$, of an algorithm as
\begin{equation}
\label{eq:tau_def}
\log \tau \equiv \lim_{t\to \infty}
\max_{\| \vz_0 \|  \le 1 } t^{-1}
\log \| \vz^t - \vz^\star \| ,
\end{equation}
where $t = 0,1,\dotsc$ is the iteration time and
$\vz^\star$ is a minimizer of the objective
function $f(\vz)$ such that $\vz^t\rightarrow \vz^*$. From \eqref{eq:tau_def}, the convergence time until $\| \vz^t - \vz^\star \|  \leq \epsilon$, for large $t$, is
\begin{equation}\label{eq:conv_time_from_rate_bound}
t = \mathcal{O}\left( \log(1/\epsilon) / \log(1/\tau) \right).
\end{equation}
If $f$ is convex and has Lipchitz gradients, then this is also an upper bound on the convergence time for $f(\vz^t) - f(\vz^\star) \leq \epsilon$.
Unless stated otherwise, $\|\cdot\|$ denotes Euclidean norm for vectors, or the Frobenius norm for matrices.

\section{Distributed Algorithms}
\label{sec:distributed_algos}

For simplicity we only discuss synchronous algorithms---many of them have asynchronously analogs which typically come with a more stringent set of assumptions to guarantee convergence. %
The algorithms we consider assume there is an agent at each node of $\G$ able to perform computations and exchange messages with its neighbours via the edges of $\G$.
To simplify the exposition, and when clear from the context, we drop the word ``agent,'' e.g., instead of  ``agent $i$ sends message to agent $j$'' we  say ``$i$ sends message to $j$.''

\paragraph{Distributed gradient descent (GD)}
We illustrate a simple distributed implementation of GD for problem \eqref{eq:pair_interaction_problem}. This 
has exactly the same iterates as its
non-distributed counterpart, $\vz^{t+1} = \vz^{t} - \alpha \nabla f\big(\vz^t\big)$, where $\alpha > 0$ is the step size. However, in such a distributed version each $e\in \E$ is responsible for computing $\nabla f_e$, while each $i\in \V$ is responsible for collecting and adding the components of all local gradients associated with the variable $z_i$ from neighboring edges. %
We summarize a distributed implementation of GD in Algorithm~\ref{alg:GD}. Note that
 $\vn_{\bm e} \equiv \vn_{N_e}$ has the same  dimension as $\vz_{\bm N_e}$ (see \eqref{eq:zNe}).
In our description, we have agents at each node/edge, but with a few redundant variables/computations it is possible to reformulate the algorithm with agents only on the nodes, or only on the edges. In Algorithm \ref{alg:GD}, the variable $z_i$ is associated to $i\in \V$ and the variable  ${\bm g}_{e}$ is associated to  $e\in \E$.
\begin{algorithm}
\SetAlgoLined
 Choose $\alpha > 0$ and initial iterate $\vz^0 \in \mathbb{R}^{|\V|}$; Set $t=0$;\\
 \While{convergence condition not met}{
  Each $i \in \V$ sends message $z^t_i$ to all edges $e \in N_i$;\\
  Each $e \in {\E}$ receives messages
  $\vz^t_{\bm e} \equiv z^t_{\bm N_e}$,
  updates $\bm{g}_e^{t+1} = \nabla f_e\big(\vz_{\bm e}^t\big)$, and sends messages $\{g^{t+1}_{e,i}\}$ to all nodes $i \in N_e$;\\
  Each $i \in \V$ receives messages $\{g^{t+1}_{e,i}\}$ from all  $e\in N_i$, and updates $z_i^{t+1} = z_i^{t} - \alpha \mbox{$\sum_{e \in N_i}$} g_{e,i}^{t+1}$;\\
  Increment $t$;\\
 }
 Read solution from $\vz^t$.
 \caption{
 Distributed GD.}\label{alg:GD}
 \end{algorithm}

Existing convergence results for distributed GD show that it is suboptimal in its dependency on $\G$ compared with the optimal algorithms discussed in Section~\ref{sec:introduction}.
For example, \cite{nedic2018network} studies a distributed subgradient method for problem \eqref{eq:consensus}, and a convex function $f$, and provides a tuning rule that leads to the following ergodic upper bound on the convergence time:
%
\begin{equation}
t = \mathcal{O}\left( (\bar{L}^2 / \epsilon^2) \left(((1 + \hat{\omega}_{\delta})/\hat{\omega}_{\delta})\right)  \right)
\end{equation}
where $\epsilon \geq f(\bar{\vz}^t) - \min_{\vz} f(\vz)$, $\bar{\vz}^t = \frac{1}{t} \sum^t_{s=1} \vz^s$,  
and $\bar{L}$ is an upper
bound on the norm of the subgradients of all 
$f_i$'s.

\paragraph{Distributed ADMM}
For distributed ADMM 
we follow \cite{Bento1,Bento2} but we further introduce a relaxation parameter $\gamma \in (0,2)$ to improve convergence, besides the usual penalty parameter
$\rho > 0$ (see \cite{Boyd:2010,FrancaBento}). Our
distributed ADMM implementation for problem \eqref{eq:pair_interaction_problem} is shown in Algorithm~\ref{alg:ADMM}. We again have agents on both nodes/edges, but it is possible to modify Algorithm~\ref{alg:ADMM} to have agents only on nodes or only on edges. The variables $x_{e,i}$ and $u_{e,i}$ (for all $i\in N_e$) sit on $e\in\E$ and $z_i$ sits on $i \in \V$.
Most convergence results for ADMM-type algorithms are not explicit in terms of properties of $\G$. Very few papers consider its behaviour with respect to the topology of $\G$, and
only loose upper bounds are known;
e.g.,
 \cite{makhdoumi2014broadcast,makhdoumi2017convergence} study 
 \eqref{eq:consensus} for smooth strongly convex functions and a tuning rule is provided leading 
 to the following 
 bound on the number of iterations required to achieve $\|\vz^t - \vz^*\| \leq \epsilon$:
 \begin{equation}\label{eq:ADMM_conv_bound_Ma}
 t = \mathcal{O}\left(\log\left({1}/{\epsilon}\right) ({\sqrt{\kappa_f}}/{\omega_L})({d^2_{\max}}/{d_{\min}}) \right)
 \end{equation}
 where $\kappa_f$ is the condition number of the objective function (i.e., the ratio of its maximum curvature and its minimum curvature).
 We know from our analysis of the canonical problem that this is suboptimal in its dependency with $\omega_L$---for a properly tuned ADMM---which should rather be  of $\mathcal{O}(1/\sqrt{\omega_L})$. We also know that there are algorithms \cite{seaman2017optimal} which, for the same class of functions,  have an $\mathcal{O}(1/\sqrt{\omega_L})$ dependency.
\begin{algorithm}
\SetAlgoLined
 Choose $\rho, \gamma > 0$, initial iterates $\vz^0,\vz^1 \in \mathbb{R}^{|\V|}$, and $\vu^0,\vx^1 \in \mathbb{R}^{|\hat{\E}|}$; Set $t = 1$;\\
 \While{convergence condition not met}{
  Each $i \in \V$ sends message $z^t_{i}$ to all $e \in N_i$;\\
  Each $e \in {\E}$ receives messages $z^t_{N_e}$, updates $u_{e,i}^{t} = u_{e,i}^{t-1} + \gamma x_{e,i}^{t} - z_{i}^{t} +
(1-\gamma) z_{i}^{t-1}$,
  updates $\vx_{\bm e}^{t+1}  = \mbox{$\argmin_{\vx_{\bm e}}$}
 f_e\left( \vx_{\bm e} \right) +
\tfrac{\rho}{2} \mbox{$\sum_{i \in N_e}$} \left( x_{e,i} - (z^t_i - u_{e,i}^t )  \right)^2$, and sends messages
$m_{e,i}^{t+1} = \gamma x_{e,i}^{t+1} + u_{e,i}^t$ to all nodes $i \in N_e$;\\
  Each $i \in \V$ receives messages $\{m^{t+1}_{e,i}\}$ from all edges $e\in N_i$, and updates $z_{i}^{t+1} = (1-\gamma)z_{i}^{t} + \tfrac{1}{|N_i|}\mbox{$\sum_{e \in N_i}$} m^{t+1}_{e,i}$ ;\\
  Increment $t$;\\
 }
 Read solution from $\vz^t$.
 \caption{
 Distributed and relaxed ADMM.}\label{alg:ADMM}
\end{algorithm}

\paragraph{Primal-dual method of multipliers (PDMM)}
An algorithm that is related to, but different from, ADMM is PDMM \cite{zhang2017distributed},   described in Algorithm \ref{alg:PDMM} in the case of problem \eqref{eq:consensus}. In PDMM, the variables $\vx_i$ and $\{\vu_{i,j}\}_{j\in N_i}$ are associated to $i \in \V$.
Under the assumptions that $f$ is strongly convex and with Lipschitz continuous gradients,
based on \cite{sherson2018derivation} one can extract, after some work, an upper bound on
the convergence time of PDMM to reach an $\epsilon$-close state to the minimizer.
%
Specifically, there is a tuning rule such that
\begin{equation}
t = \log\left( {1}/{\epsilon} \right) ({c}/{\hat{\omega}_{\delta}}) \left( 1 + \mathcal{O}\left( \hat{\omega}_{\delta})\right)\right),
\end{equation}
where $c>0$ is a constant that depends on $\kappa$, 
$d_{\max}$, and $d_{\min}$.
\begin{algorithm}
\SetAlgoLined
 Choose $\rho, \alpha > 0$ and initial iterates $\vx^0_i,\vu^0_{i,j} \in \mathbb{R}^{p}$, for all $i,j\in \V$; Set $t=0$;\\
 \While{convergence condition not met}{
  Each $i \in \V$ receives $\vu^t_{j,i}$  from its neighbour $j\in N_i$ and computes $\vn_i = \frac{1}{|N_i|} \sum_{j \in N_i}( \vx^t_j - \vu^t_{j,i})$;\\
  Each $i\in \V$ updates $\vx^{t+1}_i = \arg\min_{\vx} f_i(\vx) + \frac{\rho |N_i|}{2} \| \vx - \vn_i\|^2$;\\
  Each $i \in \V$ updates $\vu^{t+1}_{i,j} = (1-\alpha)\vu^{t}_{i,j} - \alpha(\vu^t_{j,i} + \vx^{t+1}_i - \vx^{t}_j)$, and sends and $\vu^{t+1}_{i,j}$ to all $j\in N_i$;\\
  Increment $t$;\\
 }
 Read solution from $\vx^t$.
 \caption{
 Primal-Dual Method of Multipliers.}\label{alg:PDMM}
\end{algorithm}

\paragraph{Multi-step dual accelerated methods (MSDA)}
Recent work \cite{seaman2017optimal,scaman2018optimal} provides algorithms whose convergence depends optimally on the topology of $\G$ to leading order. In Algorithm~\ref{alg:MSDA}, we describe the proposal of \cite{scaman2018optimal} for a nonsmooth convex problem \eqref{eq:consensus}. As originally proposed, at each iteration, the algorithm approximates the calculation of a proximal map for each $f_i$ using gradient information and $M<\infty$ inner iterations. However, we describe the algorithm with this map computed exactly (i.e., $M = \infty$). For the input matrix $W$ they recommend using the Laplacian $L$. The authors provide a tuning rule for all of the other constants based on the spectrum of $W$, the Lipschitz constant $\bar{L}_i$ for the gradient of $f_i$, and an upper bound $R$ on the norm of the minimizer.
Variables $\theta_i$ and $y_i$ are associated to  $i\in \V$.
The algorithm depends on a gossip subroutine that averages gradient-related information in a distributed manner;
for details on this  subroutine, including its dependency on
the parameters $c_i$ and $W$, we refer to \cite{seaman2017optimal}.
%
\begin{algorithm}
\SetAlgoLined
 Choose $W \in \mathbb{R}^{n\times n}$, $K \in \mathbb{N}$, $c_1, c_2,c_3, \eta,\sigma > 0 $, and set the initial iterates ${\bm \theta}_i^0={\bm \theta}^1_i={\bm y}_i^0 = 0 \in \mathbb{R}^{p}$ for all $i \in \{1,\dots,n\}$; Set $t=1$;\\
 \While{convergence condition not met}{
  Each $i \in \V$ computes $2{\bm \theta}^t_i - {\bm \theta}^{t-1}_i$;\\
  Starting at the above values, all agents perform $K$ rounds of gossip exchanges among themselves, at the end of which each $i\in \V$ holds a vector ${\bm g}_i \in \mathbb{R}^p$;\\
  Each $i \in \V$ updates ${\bm y}^{t+1}_i = {\bm y}^{t}_i - \sigma {\bm g}_i$;\\
  Each $i\in \V$ updates ${\bm \theta}^{t+1}_i = \argmin_{\bm \theta} f_i({\bm \theta}) + \frac{1}{2\eta} \|{\bm \theta} - (\eta {\bm y}^{t+1}_i + {\bm \theta}^t_i)\|^2$;\\
  Increment $t$;\\
 }
 Read solution from $\hat{\bm \theta}^t = \frac{1}{t \; n} \sum^n_{i = 1}\sum^t_{s =  1}{\bm \theta}^s_i$.
 \caption{
 Multi-Step Dual Accelerated method. \label{alg:MSDA} }
\end{algorithm}
\setlength{\textfloatsep}{4pt}

The tuning rule for MSDA provided in  \cite{scaman2018optimal} leads to
%
\begin{equation}\label{eq:fb_18_alg_2_convergence}
t = \mathcal{O} \left( ({R \bar{L}}/{\epsilon \sqrt{ \omega_{n} }}) + ({R \bar{L}}/{\epsilon})^2\right),
\end{equation}
where $\epsilon \geq f(\hat{\bm \theta}^t) - \min_{\bm \theta} f(\bm \theta)$, 
$\bar{L} = \|\{\bar{L}_i\}\|_2$, $\{\bar{L}_i\}$ is the vector of all $\bar{L}_i$'s, and $R$ and $\bar{L}_i$ have been explained above.
For the second method in \cite{seaman2017optimal}, under
smoothness and strong convexity assumptions, to get an $\epsilon$-close minimizer of \eqref{eq:consensus} requires
\begin{equation}\label{eq:fb_2017_smooth_strongly_convex}
t= \mathcal{O}\left( \log\left(1/\epsilon\right) \sqrt{\kappa_{\ell}}\left( 1 + 1/\sqrt{\omega_{n}}\right)  \right)
\end{equation}
iterations, where $\kappa_{{\ell}}$ is a local condition number for the $f_i$'s.

\paragraph{xFILTER} \cite{XFilter} considers finding the best convergence rate of any distributed algorithm when solving a smooth but nonconvex problem over $\G$ of the form \eqref{eq:consensus}.
They obtain a lower bound on the number of iterations required to be $\epsilon$-close to the minimum:
\begin{equation}
t = {\Omega}\left( ({\bar{L}}/{\epsilon}) ({1}/{\sqrt{\omega_L}}) \right),
\end{equation}
where $\bar{L}$ is a Lipschitz constant for the gradient of $f$. They propose an algorithm, xFILTER, described in Algorithm~\ref{alg:xFILTER}, whose convergence time matches this lower bound up to a polylog factor.
The algorithm has tuning parameters $Q \in \mathbb{N}$, $\Sigma \in \mathbb{R}^{|\E|\times|\E|}$ and $\Upsilon \in \mathbb{R}^{|\V|\times|\V|}$, both positive definite diagonal matrices with specified values; see details in \cite{XFilter}. The diagonal elements of $\Sigma$ are indexed by edges, and the diagonal elements of $\Upsilon$ are indexed by nodes.
Similarly to MSDA, xFILTER depends on a subroutine for distributed averaging that is based on Chebyshev filtering; we omit such details for the sake of space.
The variables $\vx_i, \tilde{\vx}_i, {\bm y}_i$ are associated to node $i \in \G$.
\begin{algorithm}
\SetAlgoLined
 Choose $Q \in \mathbb{N}$, $\Sigma \in \mathbb{R}^{|\E|\times|\E|}$ and $\Upsilon \in \mathbb{R}^{|\V|\times|\V|}$, both positive-definite diagonal matrices, for each $i\in \V$ set $\vx^{0}_i={\bm 0} \in\mathbb{R}^p$, ${\bm y}^0_i =-\Upsilon^{-2}_i \nabla f_i(\vx^{0}_i) \in\mathbb{R}^p$, and $\tilde{\vx}^0_i = \vx^0_i  - \Upsilon^{-2}_i \nabla f_i(\vx^0_i) \in\mathbb{R}^p$;
 Compute $R$ \cite[Eq. (74)]{XFilter}; \\
 \While{convergence condition not met}{
    Agents perform $Q$ rounds of gossip exchanges among themselves, using $\vx^t_i$ and ${\bm y}_i^t$ as starting values for node $i\in \V$, at the end of which each node $i$ has computed $\vx^{t+1}_i$;\\
    Each node $i\in\V$ computes $\tilde{\vx}^{t+1}_i = \vx^{t+1}_i - \Upsilon^{-2} \nabla f_i(\vx^{t+1}_i)$; \\
    Each node $i\in \V$ gets variables $\{\vx_j\}_{j \in N_i}$ from its neighbours and computes ${\bm y}^{t+1}_i = {\bm y}^{t}_i + (\tilde{\vx}^{t+1}_i - \tilde{\vx}^{t}_i) + \sum_{j \in N_i} \frac{\Sigma^2_{(i,j)}}{\Upsilon^2_{i}}(\vx^t_i - \vx^t_j)$;
    Increment $t$;
}
 Read solution from $\tilde{\bm x}^{t+1}$.
 \caption{
 xFILTER.
 \label{alg:xFILTER}
 }
\end{algorithm}
\setlength{\textfloatsep}{4pt}

\section{The Canonical Problem}
\label{sec:canonical}

A common analysis of an algorithm 
applied to the general problem \eqref{eq:consensus} produces relationships between its convergence rate and the spectrum of $\G$ that often hide details (e.g., constants) that are of practical importance. On the other hand, a critical step of many general algorithms is a distributed averaging step which can be analyzed in detail.
%
%
%
We now focus on such a  distributed averaging procedure via \eqref{eq:gen_quad} and ADMM.
%
%
Although the solution to \eqref{eq:gen_quad} is obvious, $\vz^* = c {\bm 1}$, characterizing the speed at which a distributed algorithm finds such a solution is nontrivial.
In solving \eqref{eq:gen_quad}, any algorithm of Section~\ref{sec:distributed_algos}---and actually any first-order algorithm---reduces to $\vx^{t+1} = T\vx^{t}$ for a vector $\vx$ of dimension equal or higher than $|\V|$, and a matrix $T$ that depends on the topology of $\G$ and on specifics of the algorithm, including its hyperparameters.
Each estimated solution, $\vz^t$, is a {linear} function of $\vx^t$, i.e.,
$\vz^t = P \vx^t$ for some $P$.
Since the solution to \eqref{eq:gen_quad}  is $\vz^* = c \bm{1}$, one then concludes:
%
\begin{equation}\label{eq:averaging_formula_for_T_x0}
c = (1/n){\bm 1}^T P T^{\infty} \vx^0.
\end{equation}
In other words, all these algorithms are solving a distributed (weighted) averaging problem over $\G$. For many algorithms, such as ADMM and GD, this average is unweighted.

%

For algorithms based on the formulation \eqref{eq:consensus}, the natural decomposition of \eqref{eq:gen_quad}, i.e., one $f_e$ per edge $e\in \E$, leads to message-passing algorithms not over $\G$ but over the \emph{line graph} $\G' \equiv (\V',\E')$ of $\G$, where $\V' \equiv \E$, and $(e_1,e_2) \in \E'$ if and only if $e_1 = (i,k) \in \E$ and $e_2 = (j,k) \in \E$, for some common node $k\in \V$.
%
%
The ADMM algorithm described in Section~\ref{sec:distributed_algos}, which we will analyze in detail in this section, is applied to the formulation \eqref{eq:pair_interaction_problem}.
For most algorithms (cf. Section~\ref{sec:distributed_algos}), the convergence rate depends on spectral properties of $\G$ and $\G'$ which control the mixing time of ``random walks.'' To leading order on the size of a graph, this mixing time is the same regardless whether the walk is done along the nodes or the edges of $\G$ (nodes of $\G'$). Therefore, to leading order, using $\G$ or $\G'$ yield the same results regarding the dependency of the convergence rate with the topology of the underlying graph. 
%


Deriving distributed averaging algorithms from \eqref{eq:gen_quad} is different than from  
\eqref{eq:alternative_formulation_for_avg_concensus};
%
%
the values to be averaged in \eqref{eq:alternative_formulation_for_avg_concensus} are  explicitly encoded in the problem, while for \eqref{eq:gen_quad} they come from the initialization, as one can see from \eqref{eq:averaging_formula_for_T_x0}.
Also, GD can be directly applied to \eqref{eq:gen_quad} since it is unconstrained.
We note that applying ADMM to \eqref{eq:alternative_formulation_for_avg_concensus} has been considered  \cite{ghadimi2015optimal,Ghadimi2014averaging}, and we will comment on  differences and similarities to our  approach  after Table \ref{table:summary} and in Table \ref{table:example}.

Solving \eqref{eq:gen_quad} is equivalent to   solving $L \vz = 0$, for which we could consider methods such as Conjugate Gradient (CG). Notice, however, that our goal is not to find any solution to  $L \vz = 0$, but a solution that corresponds to averaging information and which can be computed in a distributed way. Classically, CG and related methods have not been used in the distributed setting that we are interested---although recent  work on distributed CG exists \cite{xu2016distributed}. Classically, CG also requires $L$ to be positive definite \cite{hackbusch1994iterative}, in which case it has a convergence time $T = \mathcal{O}(1/\sqrt{\kappa})$, where $\kappa \equiv \lambda_{\max}(L)/\lambda_{\min}(L)$---although there are ways to circumvented this \cite{dostal2018conjugate}. In fact, a direct implementation of classical CG to  \eqref{eq:gen_quad} is numerically unstable with respect to its convergence to the average of the initial state, as required by an averaging gossip algorithm. We focus on ADMM and other general methods because they can solve \eqref{eq:gen_quad}, $Lz=0$, and general problems via message-passage algorithms, contrary to CG which is specific to linear systems only.

%



\subsection{ADMM, GD, and Spectral Graph Theory}
\label{sec:random_walks}

In order to apply ADMM and GD to problem
\eqref{eq:gen_quad}, let us
define $\vx = S \vz$ with  $S$ given by \eqref{eq:Sdef};
thus $\vx$ is indexed by the extended edges $\hat{\E}$ and each component $\vx_{e,i}$ contains copies of the $z_i$ incident on the extended edge $(e,i)$.
We introduce
the matrix
$Q_e = \left( \begin{smallmatrix} +1 & -1 \\ -1 & +1 \end{smallmatrix}\right)$
to form  $f_e(\vx_{\bm e})= \frac{1}{2} \vx_{\bm e}^\top Q_e \, \vx_{\bm e}$.
Note that
$f(\vx) = \frac{1}{2} \vx^\top Q \, \vx$ where
$Q = \diag(\dotsc, Q_e,\dotsc)$ is block diagonal.
Replacing this into Algorithm~\ref{alg:ADMM}, we can write its
updates concisely in  matrix form (see \cite{FrancaBentoMarkov,FrancaBentoTopology2} for details):
\begin{equation}
\label{eq:TA}
\vn^{t+1} = T_A \vn ^t, \qquad T_A \equiv I - \gamma (A + B -2BA),
\end{equation}
where $A \equiv (I + \rho^{-1} Q)^{-1}$ and $B \equiv S \big(S^T S\big)^{-1} S^T$. From $\vn^t$ we can read primal variables using $\vx^t = A \vn^t$ \cite{FrancaBentoMarkov,FrancaBentoTopology2}.
Note that that $S^T S = \D$. Moreover, $B$ is a symmetric orthogonal
projector, i.e., $B^2 = B = B^T$. The transition matrix $T_A$ is
\emph{not symmetric}, which renders
the analysis considerably more challenging than that of GD; e.g., it is unclear if $T_A$ is diagonalizable, and its eigenvalues can be complex.
%

The situation is considerably simpler for GD.
Replacing the above quantities into Algorithm~\ref{alg:GD} one
obtains a linear system:
\begin{equation}
\label{eq:TG}
\vz^{t+1} = T_G \vz^t, \qquad T_G \equiv I -\alpha L,
\end{equation}
where $L$ is the Laplacian of $\G$ (see Section~\ref{sec:notation}).
Since $T_G$ is real symmetric, it is diagonalizable and  has real eigenvalues.
Moreover, if
$\alpha \le 1/\lambda_1(L)$
then $T_G$ is also positive semidefinite and its
eigenvalues are positive. The eigenvalues of $T_G$ are
\begin{equation}
\label{eq:eigenTG}
\lambda(T_G) = 1 - \alpha \lambda(L),
\end{equation}
which is completely specified once the graph $\G$ is given.

The characterization of the eigenvalues of $T_G$, or $T_A$, directly lead to a convergence rate 
on these algorithms via a simple adaptation of a standard result from Markov chains theory.
%

%
\begin{theorem}[see \cite{MarkovChainsPeres}]\label{th:conv_markov}
Consider
$\bm{\xi}^{t+1} = T \bm{\xi}^{t}$, and let $\bm{\xi}^\star$ be a stationary
state.
If the spectral radius $\rho(T)=1$
and is attained through the eigenvalue $\lambda_1(T)=1$ with
multiplicity one,
then $\bm{\xi}^t = T^t\bm{\xi}^0$ converges to $\bm{\xi}^\star $
as $t\to\infty$ at a rate
$\| \bm{\xi}^t - \bm{\xi}^\star\| = \Theta(|\lambda_2|^t)$, where $\lambda_2 =
\lambda_2(T)$ is the second largest eigenvalue of $T$ in absolute value
(the largest is $\lambda_1(T) = 1$).
\end{theorem}
Crucial in applying Theorem \ref{th:conv_markov} for the canonical problem is the fact that  $\rho(T_A) = \rho(T_G) = 1$ and the unit eigenvalue has multiplicity one when the graph $\G$ is connected (Assumption~\ref{thm:assumption}).
One way to see this is through a detailed analysis of the spectrum of $T_A$  \cite{FrancaBentoTopology,FrancaBentoTopology2}, summarized in Lemma~\ref{thm:most_complex_eigen}; a similar analysis for $T_G$ is straightforward from \eqref{eq:eigenTG}.
An alternative argument, which can be applied to both GD and ADMM, is as follows. Let us consider ADMM. Since the canonical problem is convex, we known \cite{Boyd:2010} that the method converges in the sense that $f(\vz^t) \rightarrow \min_{\vz}f(\vz) = 0$. Note that $f(\vz) = \tfrac{1}{2}\vz^T L \vz$, where $L$ is the Laplacian of a simple connected graph, which has only one zero eigenvalue, with eigenvector the all-ones vector, ${\bm 1}$, and all its other eigenvalues are positive.  Thus, if we write $\vz^t$ as the sum of two vectors, one in the span of ${\bm 1}$ and the other in its orthogonal complement, the magnitude of the second vector goes to zero with $t$. For ADMM, since $\vn^t = A^{-1} S \vz^t$, and since ${\bm 1}$ is an eigenvector of $A$ and $S$, we have that $\vn^t$ can also be written as $\vn^t = c^t {\bm 1}  + {\bm r}^t$, where the vector ${\bm r}^t$ goes to zero with $t$, and $c^t$ is a scalar. We can thus write $(T_A)^t \vn^0 = c^t {\bm 1}  + {\bm r}^t$. Since ${\bm 1}$ is an eigenvector of $T_A$ with eigenvalue 1 (which can be checked using the definition of $T_A$ in \eqref{eq:TA}), we can make a change of basis in the last equation (using the eigenvectors of $T_A$ as a basis) and use the fact that $\vn^0$ is arbitrary to conclude that $c^{t}$ converges with $t$, and moreover,  eigenvectors that are not ${\bm 1}$ have eigenvalues with absolute value strictly smaller than 1.

Combining Theorem~\ref{th:conv_markov} with \eqref{eq:tau_def} yields
$\log |\lambda_2| = (\log \| \bm{\xi}^t - \bm{\xi}^\star \|)/t +  (\log C)/t$. Letting $t \to \infty$ we
get a convergence rate $\tau = |\lambda_2|$, controlled by a single eigenvalue of $T_G$, or $T_A$.


We now show how to compute the spectrum
of $T_A$ from the spectrum of the random walk transition matrix $\W$ of $\G$, which
has eigenvalues $\omega \equiv \lambda(\W) \in [-1,1]$. The precise analysis  \cite{FrancaBentoTopology2,FrancaBentoTopology} requires some technical precision. Here, we limit ourselves to stating the main results and providing  intuition.
\begin{lemma}[see \cite{FrancaBentoTopology,FrancaBentoTopology2}]
\label{thm:most_complex_eigen}
For each eigenvalue $\omega \equiv \lambda(\W) \in (-1,1)$ the matrix
$T_A$ has a pair of eigenvalues
\begin{equation}\label{eq:lambdaTA}
  \lambda^{\pm}(T_A) = \left(1-\dfrac{\gamma}{2}\right)
  + \dfrac{\gamma}{2+\rho}\left( \omega \pm i \sqrt{1-\dfrac{\rho^2}{4} - \omega^2}\right).
\end{equation}
Conversely, every eigenvalue $\lambda(T_A)$ has the form \eqref{eq:lambdaTA}
for some eigenvalue $w$, 
which can also be $\pm 1$ besides the
range $(-1,1)$.
\end{lemma}

Lemma~\ref{thm:most_complex_eigen} is a fundamental result to find how $\omega$, which captures
the topology of the graph $\G$, affects the convergence and the parameters of ADMM (this is the content of
Theorem~\ref{thm:tuning_ADMM} below).
For now, let us provide some intuition behind  \eqref{eq:lambdaTA}. One can write \eqref{eq:TA} in
the form
\begin{equation}\label{eq:TA2}
T_A = \left( 1 -\dfrac{\gamma}{2}\right) + \dfrac{\gamma}{\rho+2} U, \qquad
U = \Omega + \dfrac{\rho}{2} \tilde{B},
\end{equation}
where $\Omega = \tilde{B} R$ with $R = I - Q$ and $\tilde{B} = 2B - I$.
A simple matrix algebra shows that $\Omega$ is orthogonal, and both $\tilde{B}$
and $R$ are symmetric and satisfy $\tilde{B}^2 = I$ and $R^2 = I$.
The
spectrum of an orthogonal matrix lies on the unit circle in
the complex plane. Since ${B}$ is a projector, the eigenvalues of $\tilde{B}$ are $\pm 1$.
%
When $\rho = 0$ the eigenvalues of
$T_A$ lie on a unit circle, with center at $1-\gamma/2$. As $\rho$ increases,
we have a perturbation of this circle given by the term
$\tfrac{\rho}{2} \tilde{B}$.
Exploiting the structure of these matrices, it is possible to show that the eigenvalues
of $T_A$ still lie on a circle, again centered at
$(1-\gamma/2)$, but now with radius
$\tfrac{\gamma}{2} \sqrt{(2-\rho)/(2+\rho)}$ which decreases with $\rho$.

An important eigenvalue of $T_A$ that might compete for $|\lambda_2|$, and which will in general lie outside the aforementioned circle, is the one studied in following.
\begin{lemma}[see \cite{FrancaBentoTopology2}] \label{thm:real_eigen}
The matrix $T_A$ has eigenvalue $\lambda(T_A) = 1-\gamma$ if and only
if $\G$ has at least one cycle of even length.
\end{lemma}
Notice that, for large $n$, many random graph models
have even cycles with high probability.
In Section \ref{sec:explicit_conv_rage}, Theorem~\ref{thm:most_complex_eigen} and Lemma \ref{thm:real_eigen} are used to obtain exact rates for ADMM.


\subsection{Markov Chain Lifting}

We now show a surprising connection between ADMM and \emph{lifted
Markov chains} \cite{FrancaBentoMarkov}; we note that a connection between GD and Markov chains  is well-known in the literature \cite{BertsekasBook}.

%
Let $\MC$ and $\hMC$ be two finite Markov chains with states
$\V$ and $\hV$, respectively. Assume that $|\V| < |\hV|$.
The respective transition probability matrices are denoted by
$M$ and $\hM$, and their stationary distributions by $\vpi$ and $\hvpi$.
One would like to sample $\vpi$ from $\hvpi$. A set of conditions when
this is possible is known as \emph{lifting}. We say that $\hMC$ is
a lifting of $\MC$ if there exists a row stochastic matrix
$S \in \mathbb{R}^{|\hV|\times |\V|}$ such that
\begin{equation}\label{eq:lifting}
\vpi = S^T \hvpi, \qquad D_{\vpi} M = S^T D_{\hvpi} \hM S.
\end{equation}
Above, $D_{\vpi} \equiv \diag(\pi_1,\dotsc,\pi_{|\V|})$ and similarly for
$D_{\hvpi}$. What this means is that it is possible to collapse
the larger Markov chain $\hMC$ onto the smaller Markov chain $\MC$.
This is interesting
when $\MC$ is slow mixing but $\hMC$ is fast mixing so that one
can sample $\vpi$ efficiently through $\hvpi$. The mixing
time $\HH$ of $\MC$ is a measure of the time it takes to approach
stationarity, e.g., 
\begin{equation} \label{eq:mixing}
\HH \equiv \min\Big\{ t : \max_{i,\bm{p}^0} |p_i^t - \pi_i| < 1/4\Big\},
\end{equation}
where $p_i^t$ is the probability of hitting
state $i$ after $t$ steps starting with distribution $\bm{p}^0$.
For walks on a ring, a $2D$-grid, or a $D$-dimensional grid, lifted chains have been produced \cite{lovasz1999faster} whose  mixing time obeys
\begin{equation}\label{eq:mixing_sqrt}
\hHH \approx C \sqrt{\HH},
\end{equation}
where $C \in (0,1)$ is some constant that depends on $\vpi$. In many other cases, lifting achieves nontrivial speedups  \cite{Diaconis}. For some family of random walks, general lifting procedures have also been constructed  \cite{lovasz1999faster}.
However, there are formal results establishing a limit on the attainable
speedup \cite{chen1999lifting}: if $\MC$ is irreducible then
$\hHH \ge C \sqrt{ \HH } \log(1/\pi_0)$, where $\pi_0$ is the smallest stationary probability of any state, and if both $\MC$ and $\hMC$ are reversible
then $\hHH \log(1/\pi_0) \ge C \HH$.

Let us now introduce two matrices $M_G$ and $M_A$ that are closely
related to $T_G$ and $T_A$, respectively. We also introduce the associated state vectors $\vv_G$ and $\vv_A$:
%
\begin{align}
  M_G &\equiv (I - D_G)^{-1}(T_G - D_G), & \vv_G  &\equiv (I - D_G) \bm{1}, \label{eq:MG} \\
  M_A &\equiv (I - D_A)^{-1}(T_A - D_A), & \vv_A  &\equiv (I - D_A)\bm{1}.\label{eq:MA}
\end{align}
Above, $D_G \ne I$ and $D_A \ne I$ are arbitrary diagonal matrices, to be specified
shortly.
%
%
It is possible to show that $M_G$ is the transition matrix
of a Markov chain, while $M_A$ ``almost'' defines a
Markov chain, i.e., $M_A$ has all the necessary properties of a
transition matrix except that it can have few
``negative entries.''

\begin{theorem}[see \cite{FrancaBentoMarkov}]
The following holds:
\begin{enumerate}
\item For $(D_G)_{ii} < 1$ and sufficiently small $\alpha$, $M_G$ is a doubly
stochastic matrix, hence a probability transition matrix.
\item The rows of $M_A$ sum up to one. There exists a graph $\G$ such
that, for any choice
of $D_A$, $\rho$ and $\gamma$, at least one entry of $M_A$ is negative.
\item We have that $\vv_G^T M_G = \vv_G^T$ and $\vv_A^T M_A = \vv_A^T$,
so these vectors can be seen as (right) stationary distributions of
$M_G$ and $M_A$, respectively.
\end{enumerate}
\end{theorem}

Therefore, except for  potentially few negative entries, $M_A$ acts as Markov chain.
The next surprising result shows that there is a well-defined
lifting between ADMM and GD.

\begin{theorem}[see \cite{FrancaBentoMarkov}] \label{thm:lifting}
The matrix $M_A$ is a lifting of $M_G$ in the sense of \eqref{eq:lifting}.
Specifically:
\begin{equation} \label{eq:liftingAG}
\vv_G = S^T \vv_A, \qquad D_{\vv_G} M_G = S^T D_{\vv_A} M_A S,
\end{equation}
provided the following relations hold true:
\begin{equation} \label{eq:lifting_params}
\rho S^T (I - D_A) S = I - D_G, \qquad
\alpha = \gamma \rho / (\rho + 2).
\end{equation}
\end{theorem}

\begin{figure}
\centering
\includegraphics[scale=0.8]{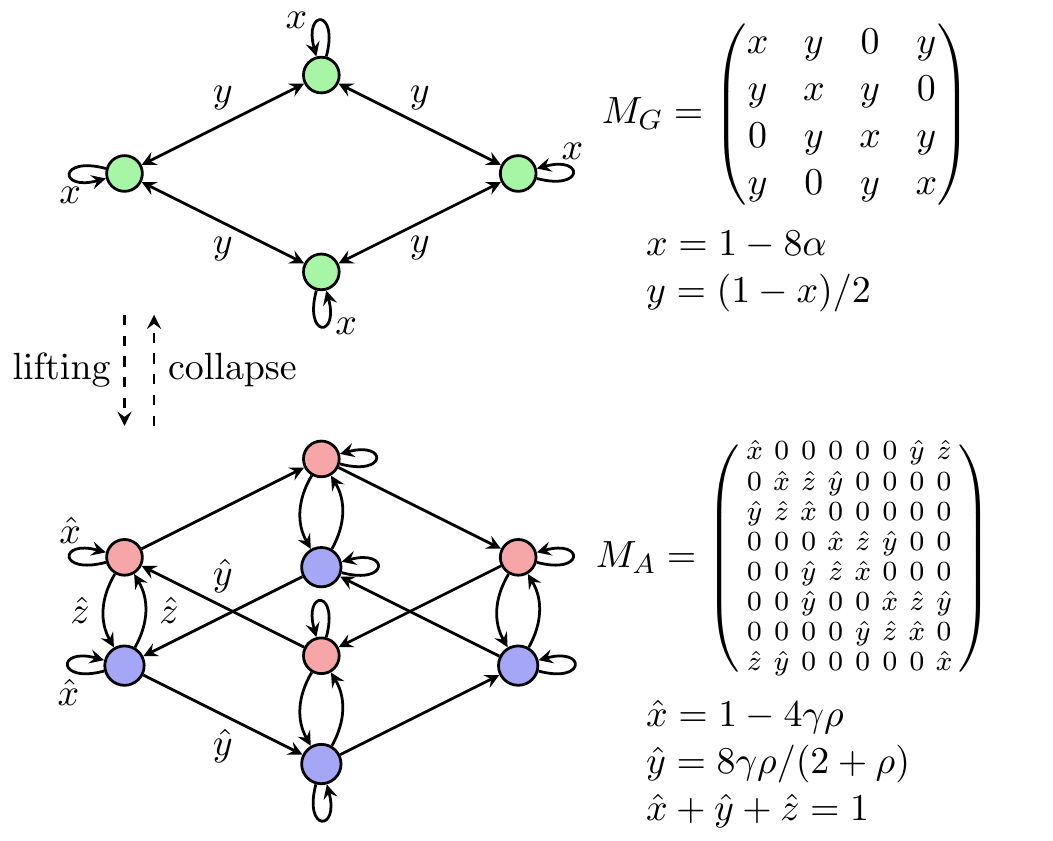}
\caption{\label{fig:ring_lifting}
Lifting between ADMM and GD for a ring graph.
We obtain  $M_A$ and $M_G$ by solving \eqref{eq:lifting_params} and it is easy
to check \eqref{eq:liftingAG}.
Here $\vv_G = \tfrac{1}{4} \bm{1}$ and
$\vv_A = \tfrac{1}{8} \bm{1}$.
Note that the lifted chain has two opposite directions
for `lower'' and ``upper''
graphs as consequence of $M_A$ being nonsymmetric.
}
\end{figure}

In \cite{FrancaBentoMarkov} this is actually proven for slightly more
general problems than \eqref{eq:gen_quad} that include weights in each term, and with
ADMM allowed to have several penalty parameters.
%
%
In Fig.~\ref{fig:ring_lifting} we show an explicit example of lifting between
GD and ADMM for a ring graph; in this case all entries of $M_A$ are actually positive
and this is a lifting in the truly Markov chain sense.
We note that the above relation between $\alpha$, $\gamma$, and $\rho$ is necessary
to establish the lifting connection, but 
we stress that this relation \emph{will not be assumed} in the following results.

Since in several cases Markov chain lifting provides a square root speedup \eqref{eq:mixing_sqrt} (see e.g. \cite{chen1999lifting,Diaconis}), and further
due to Theorem~\ref{thm:lifting}  together with numerical evidence, we conjectured
\cite{FrancaBentoMarkov} that actually a stronger inequality is attained
between ADMM in comparison to GD. The convergence rate $\tau$ is related
to the mixing time \eqref{eq:mixing} as $\HH \sim 1/(1-\tau)$. Thus,
denoting $\tau_G^\star$ and $\tau_A^\star$ the optimal convergence rates
of GD and ADMM (obtained with optimal parameters), respectively, we
conjectured that there exists
a constant $C > 0$ such that
\begin{equation}
  \label{eq:conjecture}
1-\tau_A^\star \ge C \sqrt{1-\tau_G^\star}.
\end{equation}
This is equivalent to $\HH_A \le C \sqrt{\HH_G}$. Note that this is
a  stronger condition than for lifted Markov chains where
one usually has, at best, $\hHH \ge C \sqrt{\HH}$. Moreover, this statement claims that \eqref{eq:conjecture} must always holds, even for graphs
with bottlenecks for which Markov chain lifting does not 
speedup.\footnote{
The relation \eqref{eq:conjecture} was actually established in followup
work \cite{FrancaBentoTopology2,FrancaBentoTopology} and will be described below.}

\subsection{Explicit Convergence Rate}\label{sec:explicit_conv_rage}

The lifting relation previously discussed provides a surprising mathematical connection,
suggesting that the significant speedup of ADMM in comparison to GD has  origin in the fact that ADMM can be seen as a ``higher-dimensional'' version of GD.
Now we explicitly show ADMM's optimal rate by optimizing
the eigenvalues of $T_A$. We also  compare such a rate with the optimal rate of GD. Thus, we  provide explicit formulas
for ADMM's optimal convergence and further show that  \eqref{eq:conjecture} holds true.

Note that \eqref{eq:lambdaTA} depends on the eigenvalues
$\omega \in [-1,1]$ of the transition matrix $\W$. The second
largest eigenvalue, denoted by $\omega^\star \equiv \lambda_2(\W)$,
plays an important role; this eigenvalue is related to the mixing time
of $\W$ and also to the conductance $\Phi \in [0,1]$ of $\G$
through the Cheeger bound \cite{Cheeger}:
\begin{equation} \label{eq:conductance}
  1-2 \Phi \le \omega^\star \le 1-\Phi^2/2.
\end{equation}
High conductance means fast mixing, while low conductance means
slow mixing and indicates the presence of bottlenecks. 
Based on the lifting connection, the most
interesting cases are the ones with low conductance where in the Markov
chains world even the lifted chain cannot speedup over the base chain.
However,
we conjectured in \eqref{eq:conjecture} that even in these cases ADMM
should improve over GD. Thus,  we will
later focus on the case $\omega^\star \in [0, 1)$ which means
$\Phi \le 1/2$.

According to Theorem~\ref{th:conv_markov},
to tune ADMM one needs to minimize the second largest eigenvalue of
$T_A$---in absolute value.
Such an eigenvalue comes either from \emph{(a)} the conjugate pairs in
\eqref{eq:lambdaTA} with $\omega^\star = \lambda_2(\W)$,
or from \emph{(b)} the real eigenvalue $1-\gamma$ of Lemma~\ref{thm:real_eigen} (see
\cite{FrancaBentoTopology2} for more details).
The complex eigenvalues lie on a circle in the complex plane, whose radius
shrinks as $\rho$ increases. Thus, we must increase
$\rho$ to make the radius as small as possible, which happens
when the complex eigenvalues \emph{(a)} fall on the real line.
This will determine the tuning $\rho^\star$.
Now, considering the real eigenvalue \emph{(b)}, we can fix $\gamma^\star$ by making
$|1-\gamma|$ the same size as the norm of the previous complex conjugate
eigenvalues that just felt on the real line.
These ideas lead to the following result.

\begin{theorem}[see \cite{FrancaBentoTopology2}]
\label{thm:tuning_ADMM}
Let $\W$ be the random walk transition matrix of $\G$,
 $\omega^\star$ the second largest (not in absolute value) eigenvalue of
$\W$, and $\bar{\omega}$ the smallest (not in absolute value) eigenvalue of $\W$ that is different than $-1$.
Let $\lambda_2(T_A)$ be the second largest
eigenvalue of $T_A$ in absolute value. The optimal convergence rate of ADMM is given by
\begin{equation}
\label{eq:optimal_rate_ADMM}
\tau^\star_A \equiv \min_{\gamma, \rho} |\lambda_2(T_A)|, 
\end{equation}
with parameters $\gamma^\star$ and $\rho^*$ provided in Table~\ref{table:summary}.
%
%
\end{theorem}

This result explicitly describes the behaviour of ADMM in terms
of spectral properties of the underlying graph $\G$.
We recall
that besides the Cheeger bound \eqref{eq:conductance},
$\omega^\star$ is related to the well-known \emph{spectral gap}, which
determines the algebraic connectivity of the graph. In Section \ref{sec:numerics}, we provide values of $\omega^\star$ for some common graphs, which together with Table \ref{table:summary} and \eqref{eq:conv_time_from_rate_bound} produce examples of convergence times. In a nutshell, the more connected $\G$ is, the smaller is the convergence time (for a fixed solution accuracy).


\begin{table}[h]
\caption{\label{table:summary} We summarize the optimal convergence rate
of ADMM and optimal parameters under different constraints on the topology
of the graph $\G$.
The most interesting and common case is the one in the first table.
The constant $\omega^\star$ is the second largest (not in absolute value) eigenvalue of the transition matrix
$\W$, while $\bar{\omega}$ is the smallest (not in absolute value) eigenvalue that is different than $-1$.}

\centering

\begin{tabular}{@{}l|ll@{}}
\multicolumn{3}{c}{(a) $\G$ has even lengh cycles.} \\
\toprule[1pt]
& $\omega^\star \ge 0$ & $\omega^\star < 0$ \\
\midrule[.5pt]
$\rho^\star$ & $2\sqrt{1-(\omega^\star)^2}$ & $2$ \\
$\gamma^\star$ & $4\left(3-\sqrt{\tfrac{2-\rho^\star}{2+\rho^\star}}\right)^{-1}$
&  $4/3$ \\
$\tau_A^\star$ & $\gamma^\star - 1 = 1-\sqrt{2} \sqrt{1-\omega*}+\mathcal{O}\left(1-\omega*\right)$ & $\gamma^\star - 1$ \\
\bottomrule[1pt]
\end{tabular}

\vspace{1em}

\begin{tabular}{@{}l|l l l @{}}
\multicolumn{4}{c}{(b) $\G$ has a cycle, but not with an even length.}\\
\toprule[1pt]
& $ 0 \le \omega^\star \le |\bar{\omega}|$ &
$0 \le |\bar{\omega}| < \omega^\star $ & $\omega^\star < 0$ \\
\midrule[.5pt]
$\rho^\star$ &
$2\sqrt{1-(\omega^\star)^2}$ &
$2\sqrt{1-(\omega^\star)^2}$ &
$2$ \\
$\gamma^\star$ &
$\tfrac{2(2+\rho^\star)}{2+\rho^\star-\bar{\omega}-\omega^\star +
\sqrt{\bar{\omega}^2 - (\omega^\star)^2}}$ &
$2$ &
$4/(2-\bar{\omega})$ \\
$\tau_A^\star$ &
$1-\gamma^\star\left(\tfrac{1}{2} - \tfrac{\omega^\star}{2+\rho^\star}\right)$ &
$\tfrac{2\omega^\star}{2+\rho^\star}$ &
$1-\gamma^\star\left(
\tfrac{1}{2}-\tfrac{\omega^\star}{2+\rho^\star} \right)$ \\
\bottomrule[1pt]
\end{tabular}

\vspace{1em}

\begin{tabular}{@{}l|l l l @{}}
\multicolumn{4}{c}{(c) $\G$ does not has cycles.}\\
\toprule[1pt]
& $ 0 \le \omega^\star \le |\bar{\omega}|$ &
$0 \le |\bar{\omega}| < \omega^\star $ & $\omega^\star < 0$ \\
\midrule[.5pt]
$\rho^\star$ &
$2\sqrt{1-(\omega^\star)^2}$ &
$2\sqrt{1-(\omega^\star)^2}$ &
$2\sqrt{1-\bar{\omega}\omega^\star}$ \\
$\gamma^\star$ &
$\tfrac{2(2+\rho^\star)}{2+\rho^\star-\bar{\omega}-\omega^\star +
\sqrt{\bar{\omega}^2 - (\omega^\star)^2}}$ &
$2$ &
$\tfrac{2+\rho^\star}{1-\bar{\omega}+\rho^\star/2}$ \\
$\tau_A^\star$ &
$1-\gamma^\star\left(\tfrac{1}{2} - \tfrac{\omega^\star}{2+\rho^\star}\right)$ &
$\tfrac{2\omega^\star}{2+\rho^\star}$ &
$\tfrac{\sqrt{\bar{\omega}(\bar{\omega}-\omega^\star)}}{
1-\bar{\omega}+\sqrt{1-\bar{\omega}\omega^\star}}$ \\
\bottomrule[1pt]
\end{tabular}
\end{table}

The results in Table \ref{table:summary} differ from \cite{ghadimi2015optimal,Ghadimi2014averaging},
for which the problem formulation is different, namely
\eqref{eq:gen_quad} versus
\eqref{eq:alternative_formulation_for_avg_concensus}, respectively.
By means of \eqref{eq:alternative_formulation_for_avg_concensus} there is no connection with lifted Markov chains. Moreover, the ADMM approach in \cite{ghadimi2015optimal,Ghadimi2014averaging} has many more tuning parameters, one per edge in $\G$, and
are required to satisfy constraints (see \cite[Assumption 1]{ghadimi2015optimal}); there is no explicit tuning rules for these parameters, although it is possible to do so via numerically solving an SDP. The optimal rates provided in \cite{ghadimi2015optimal,Ghadimi2014averaging} depend on these (unknown) tuned parameters (called edge-weights). These facts render such an approach less transparent than the ADMM formulation
emphasized in this paper.
The results of  \cite{ghadimi2015optimal,Ghadimi2014averaging} provide three different formulas for the convergence rate for three different cases (called C1, C2, and C3),
the selection of which depends on the relative magnitude of the spectral properties of some matrices. In contrast, for the results summarized in Table~\ref{table:summary},  the cyclic properties of $\G$ play a role in selecting the correct formula as well. Furthermore,  two of the cases (C2 and C3) are \emph{not} proven to be optimal  \cite{ghadimi2015optimal,Ghadimi2014averaging}.

Despite of these facts, it is possible to establish some commonalities. Because the formulas in \cite{ghadimi2015optimal,Ghadimi2014averaging} are not as transparent as those in Table~\ref{table:summary},  we simplify them using the recommendation of setting the edge-weights to the inverse the degree of the nodes (i.e., set  $\W_{i,(i,j)} = 1/d_i$). For simplicity, we further consider $\G$ to be a regular graph. In this case, the rates for two of the cases in \cite{ghadimi2015optimal,Ghadimi2014averaging} (C1 and C2) are $\tfrac{2\omega^\star}{2+\rho^\star}$ and $1-\gamma^\star\left(\tfrac{1}{2} - \tfrac{\omega^\star}{2+\rho^\star}\right)$, which also appear in Table \ref{table:summary}. When $-1$ is not an eigenvalue of $\W$, the rate for C3 is $\tfrac{-\bomega}{1-\bomega}$, and does not appear in Table \ref{table:summary}. Several of the above formulas, specially $\tau^*_A$ in Table \ref{table:summary} (a), have no analog in \cite{ghadimi2015optimal,Ghadimi2014averaging}. These differences are reflected  in the numerical results (Table \ref{table:example} below). 

Finally, we state an extended version of conjecture \eqref{eq:conjecture}. 

\begin{theorem}[see \cite{FrancaBentoTopology2}]\label{thm:proof_conj} Suppose that $\G$ has an even length
cycle and conductance $\Phi \le 1/2$. Let $\Delta \equiv d_{\textnormal{max}}/ d_{\textnormal{min}}$. There is a constant
$C = 1 - O(\sqrt{{\omega^*}_{\delta}})$, where ${\omega^*}_{\delta} \equiv 1-\omega^\star$, such that
\begin{equation}\label{eq:proof_conjecture}
C (1-\tau_G^\star) \le (1-\tau_A^\star)^2 \le 2 \Delta C (1-\tau_G^\star).
\end{equation}
\end{theorem}

The right hand side of \eqref{eq:proof_conjecture}
implies that the square root speedup attained by ADMM is tight.
Nevertheless, the gap becomes larger for very irregular graphs,
which have $\Delta \gg 1$, compared to regular graphs, which have $\Delta = 1$.
Interestingly, Theorem~\ref{thm:proof_conj} holds for
\emph{any} graph. This is in contrast
to lifted Markov chains which cannot speedup when the conductance
of the base graph is small.

The proof of Theorem~\ref{thm:proof_conj} is based on the inequalities  \eqref{eq:relation_among_spectral_measures}.
Expanding the convergence
rate of ADMM \eqref{eq:optimal_rate_ADMM} in terms of $1-\omega^*$ (see Table \ref{table:summary} (a), left column), we use \eqref{eq:relation_among_spectral_measures} to relate $1-\omega^*$ with the eigenvalues of the Laplacian and then with the convergence rate of GD   \eqref{eq:eigenTG}.
At the same time, the largest magnitude
eigenvalue of $T_G$ which is smaller than $1$  can be minimized by taking
$\alpha = 2 (\lambda_1(L) + \lambda_{|\V| -1}(L) )^{-1}$.

\bigskip

\paragraph*{Example}
Let us illustrate how to use the results of Table~\ref{table:summary}.
Given a graph $\G$, one must first check into which of the three categories it
fits, i.e., whether there are even length cycles or not, etc.
Then, one must compute $\omega^\star$ and $\bar{\omega}$ from the transition matrix of the
graph. From this, the optimal parameters $\rho^\star$ and $\gamma^\star$ of ADMM follow, and also
its optimal convergence rate $\tau_A^\star$. We will use this approach to verify numerically the
performance of ADMM in the next section.
As an example, we show some simple cases in Table~\ref{table:example}.
In particular, note that $\gamma^\star > 1$, in agreement with empirical knowledge
that over-relaxation improves performance. The penalty parameter is always
$1 < \rho^\star \le 2$.
The last row of Table~\ref{table:example}  illustrates how the formulation \eqref{eq:gen_quad} may lead to a faster distributed averaging compared to the more traditional formulation \eqref{eq:alternative_formulation_for_avg_concensus}.

\begin{table}[h]
\renewcommand\arraystretch{1.1}
\caption{
Examples in computing  optimal parameters and convergence rates
for ADMM based on
Table~\ref{table:summary}; the second to last row
indicates the subtable and  column being used. 
The last row indicates the optimal rates for the ADMM averaging algorithm of \cite{ghadimi2015optimal,Ghadimi2014averaging} to solve \eqref{eq:alternative_formulation_for_avg_concensus}, called the ``edge-variable formulation;'' we use the provided formulas to tune $\rho$ and $\alpha$, and we set each edge-weight parameters $\mathcal{W}_{i,(i,j)} = 1/d_i$ (recall that $d_i$ is the degree of node $i$), such that Assumption~1 of \cite{ghadimi2015optimal,Ghadimi2014averaging} hold.
}
\label{table:example}
\centering
\begin{tabular}{@{}l | l l l l@{}}
\toprule[1pt]
$\G$
 & \includegraphics[scale=.25]{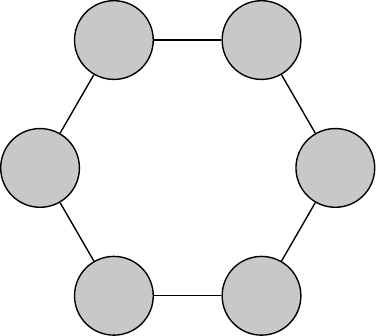}
 & \includegraphics[scale=.25]{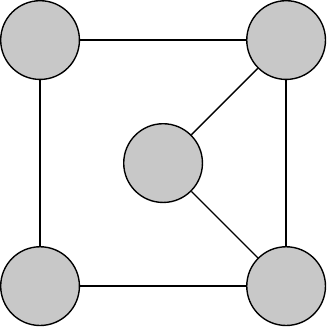}
 & \includegraphics[scale=.25]{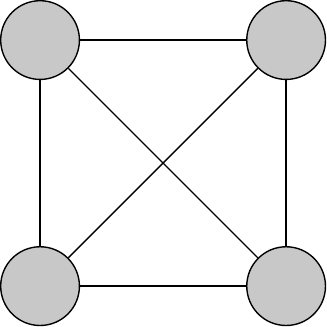}
 & \includegraphics[scale=.25]{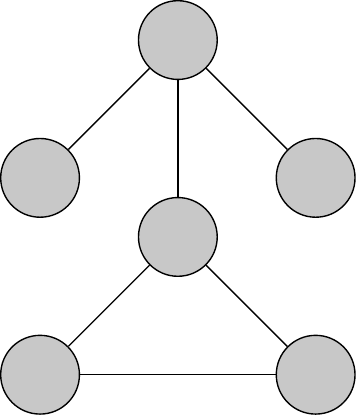} \\
\midrule[0.5pt]
$\Phi$ & $1/3$ & $1/2$ & $1$ & $1/5$ \\
$\omega^\star$ & $1/2$  & $1/3$ & $-1/3$ & $(\sqrt{97}-1)/12$ \\
$\rho^\star$ & $1.732$ & $1.886$ & $2$ & $1.351$ \\
$\gamma^\star$ & $1.464$ & $1.414$ & $4/3$ & $1.659$ \\
$\tau_A^\star$ & ${0.464}$  & $0.414$ & $1/3$ & $0.536$ \\
Table~\ref{table:summary}& (a) Col.~1  & (a) Col.~1  & (a) Col.~2  & (b) Col.~1  \\
\hline
Ref.~\cite{ghadimi2015optimal,Ghadimi2014averaging}& 0.634  & 0.594  & 1/7  & 0.761  \\
\bottomrule[1pt]
\end{tabular}
\end{table}

\section{Numerical Experiments} \label{sec:numerics}

We now consider some experiments to verify  the previous theoretical results, and also to verify the practical performance
of ADMM compared to the other state-of-the-art algorithms described in Section~\ref{sec:distributed_algos}
when solving distributed-averaging problems.
%
We consider consensus problems over a few standard graphs $\G$, namely:
\begin{itemize}
\item a graph sampled from the
Erd\" os-Renyi model with edge probability
$p = \log(n)/n$; $1/(1-\omega^\star) = \mathcal{O}(1)$;
\item a $k$-hop lattice graph with $k = \log n$. This is just
a ring with an extra edge connecting each node to another node at a distance
$k$ apart; $1/(1-\omega^\star) = \mathcal{O}(n^2/\log^2 n)$;
\item a ring graph; $1/(1-\omega^\star) = \mathcal{O}(n^2)$.
\item a periodic grid graph; $1/(1-\omega^\star) = \mathcal{O}(n^4)$.
\end{itemize}
Above, we also indicate the asymptotic value of $1/(1-\omega^\star)$ as a function
of the number of nodes $n$, which is related to the dependency of the convergence time of ADMM, and other algorithms, with respect to $\G$. We  compare  the following:
\begin{enumerate}
\item GD, given by Algorithm~\ref{alg:GD}, for solving \eqref{eq:consensus};
\item MSDA1 \cite[Algorithm~1]{seaman2017optimal}  for solving \eqref{eq:consensus};
\item MSDA2 \cite[Algorithm~2]{seaman2017optimal}  for solving \eqref{eq:consensus};
\item MSDA3, given by Algorithm~\ref{alg:MSDA}, for solving \eqref{eq:consensus};
\item PDMM, given by Algorithm~\ref{alg:PDMM}, for solving \eqref{eq:consensus};
\item ADMM1, which is a
 consensus ADMM implementation for solving \eqref{eq:consensus}, and where all agents hold a local copy of the complete $\vz$;
\item ADMM2, which is the same
as the previous algorithm, but where we do not introduce a consensus variable in the ADMM algorithm, hence it is very close to PDMM;
\item ADMM3, which is precisely
Algorithm~\ref{alg:ADMM} described in this paper, and which solves problems of the form \eqref{eq:pair_interaction_problem}, and where each agent only holds local copies of few components of $\vz$.
\item xFILTER \cite[Algorithm~\ref{alg:xFILTER}]{XFilter}.
\end{enumerate}

We consider the performance of these algorithms on two specific problems.
Some of the above algorithms have inner loops. The iteration number that we report in our experiments is the number of times that the code in the inner-most loop is executed.
We note that the code used in our experiments is available in \cite{GitLink}.

\subsection{Canonical Problem}
\begin{figure*}[t]
\centering
\includegraphics[width=0.32\textwidth,trim={16 0 40 0 },clip]{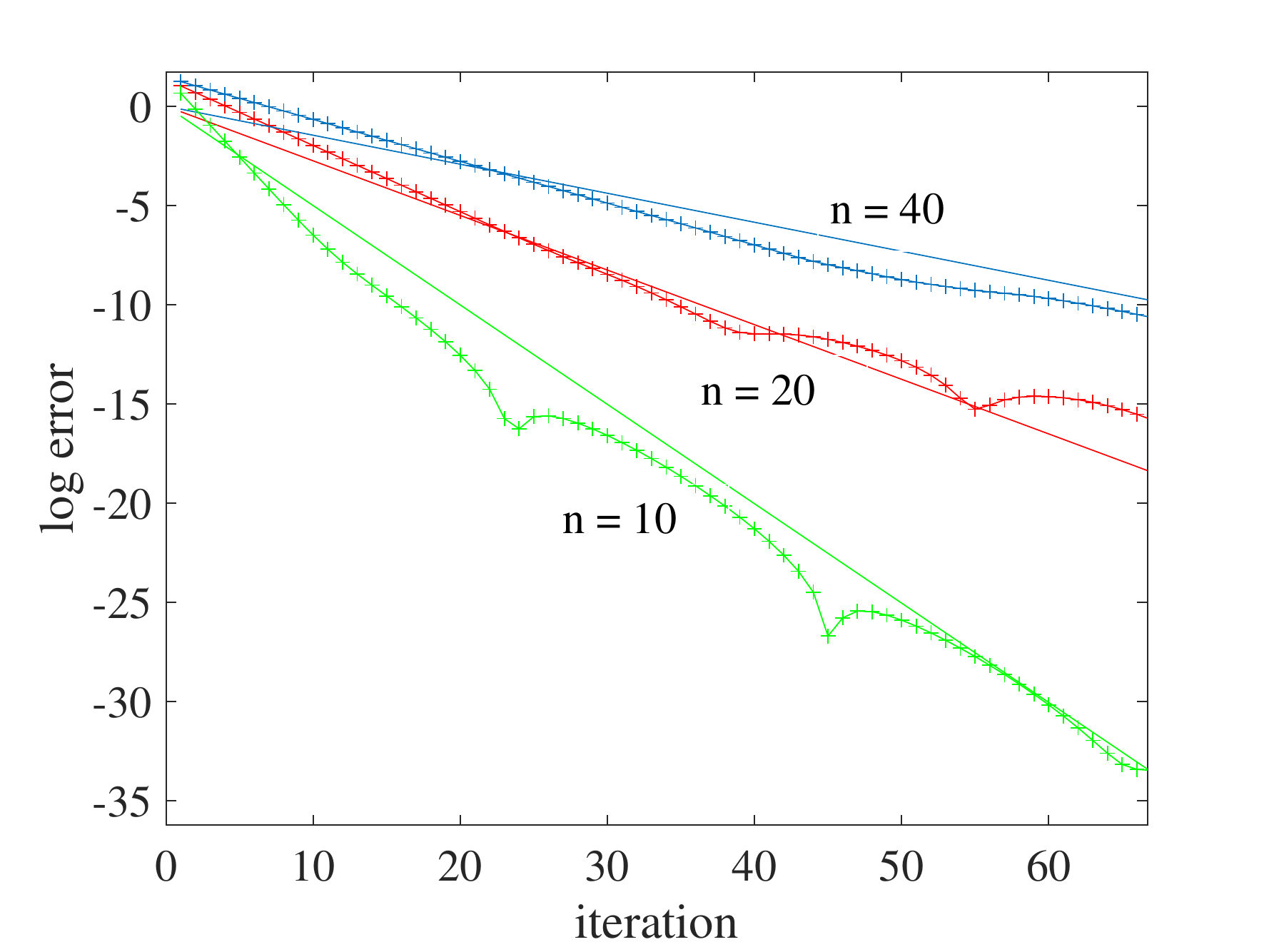}
\includegraphics[width=0.32\textwidth,trim={16 0 40 0 },clip]{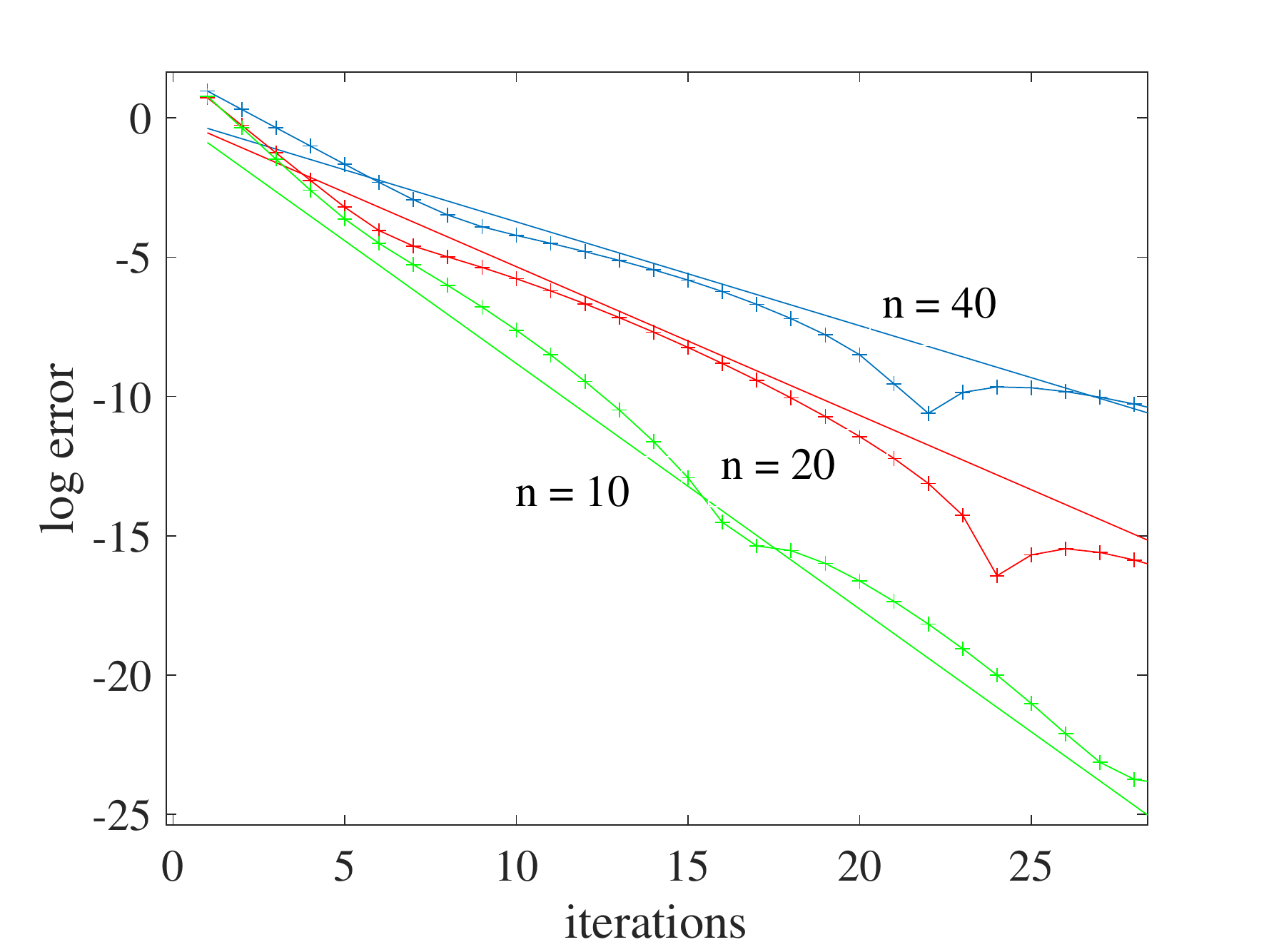}
\includegraphics[width=0.32\textwidth,trim={16 0 40 0 },clip]{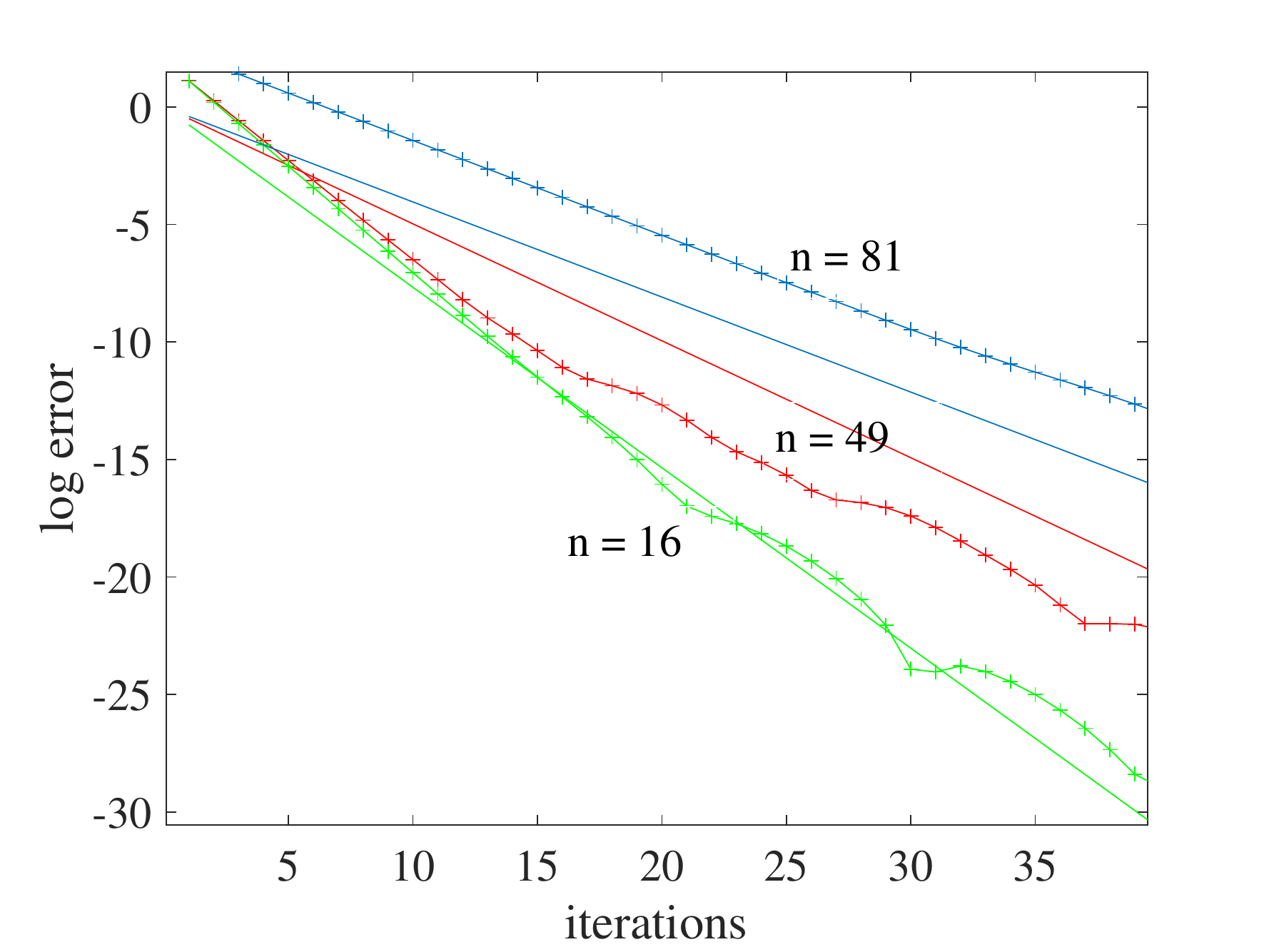}
\caption{\label{fig:canonical_test}
Comparison between theory
and an actual execution of ADMM (i.e., ADMM3) on problem \eqref{eq:gen_quad}.
The solid straight lines correspond to our theoretical predictions for the convergence rate
(see Table~\ref{table:summary}),
while the lines with markers correspond to ADMM  where the optimal parameters
were found using Bayesian optimization. These numerical results are consistent with the theoretical results since these lines have similar slopes.
\emph{Left:} ring graph. \emph{Middle:} $k$-hop lattice graph with $k = \log n$. \emph{Right:} 2D periodic grid graph.
In each panel we indicate three graphs of the same type but with different number, $n$, of nodes.
}
\end{figure*}


We compare the performance of several algorithms on problem \eqref{eq:gen_quad}, including algorithms that are designed for strongly convex problems; note that \eqref{eq:gen_quad} is only convex.
To be able to run such algorithms on the same problem, we add a  small regularization to \eqref{eq:gen_quad}, namely,  $\frac{\delta}{2}\frac{|\V|}{|\E|}\|\vz - {\bm c}\|^2$, with $\delta = 10^{-3}$ and ${\bm c} \neq {\bm 0}$.

Our first goal is to verify if the rates predicted by
Theorem~\ref{thm:tuning_ADMM}
can be achieved in practice;
it could be the case that ADMM is too sensitive to small changes in its parameters that in practice it is impossible to obtain the optimal theoretical rates.
We thus  tune ADMM3 using Bayesian optimization to find if the rates achieved in this way are actually close to the rates predicted by our theoretical formulas. 


For each run of ADMM3, we plot the error $\log \|\vz^{t} - {\bm c}\|$ versus $t$, and  compare if the slope of this curve is close to the one provided by our formulas. We see that this is indeed the case in the examples
of Fig.~\ref{fig:canonical_test}. Each panel considers the same graph but with different number
of nodes. The straight solid line is the theoretical slope, while the lines with markers
correspond to the empirical behaviour of ADMM3 tuned with Bayesian optimization.

\begin{figure}
\centering
\includegraphics[scale=.40,trim={15 0 20 0},clip]{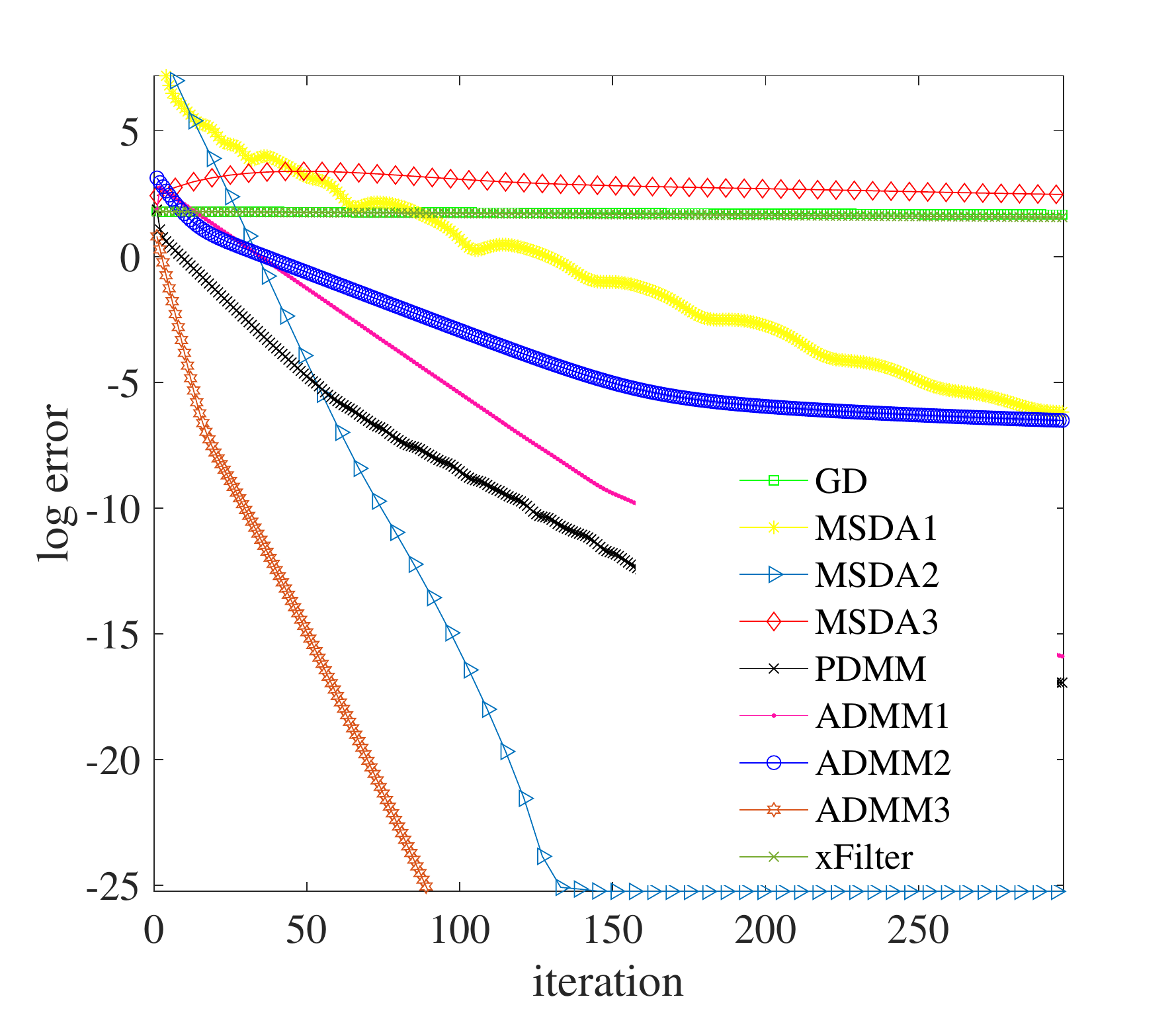}
\caption{\label{fig:all_algs_ring}
Comparison of different algorithms for problem \eqref{eq:gen_quad} and a ring graph with $n=20$ nodes. The algorithms were tuned using Bayesian optimization.
}
\end{figure}

Next, we consider problem \eqref{eq:gen_quad}---plus the small regularization term mentioned above---with a ring graph
with $n=20$ nodes, and we compare the convergence rate of different algorithms on this problem in Fig.~\ref{fig:all_algs_ring}. The  error is computed as above.
All the algorithms were tuned using Bayesian optimization.
We see that ADMM and MSDA2 outperform the other methods, however MSDA2 requires strong convexity, as opposed to  ADMM.
We note that GD, xFILTER, and MSDA3 also converge to arbitrary accuracy but their convergence time are orders of magnitude slower compared to the other algorithms.

\subsection{Sensor Localization Problem}

\begin{figure*}
\centering
\includegraphics[width=0.325\textwidth,trim={5 0 10 0},clip]{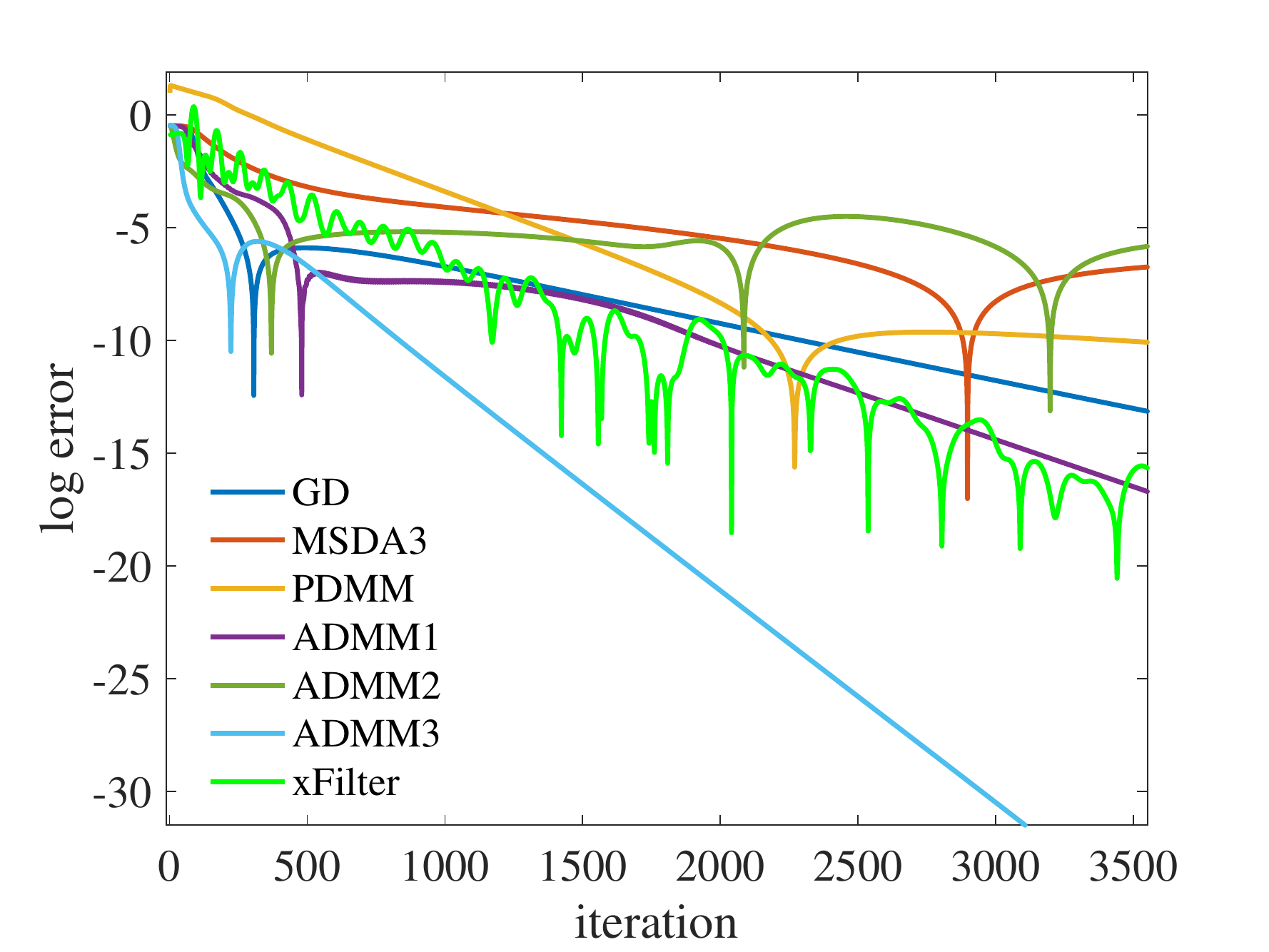}
\includegraphics[width=0.325\textwidth,trim={5 0 10 0},clip]{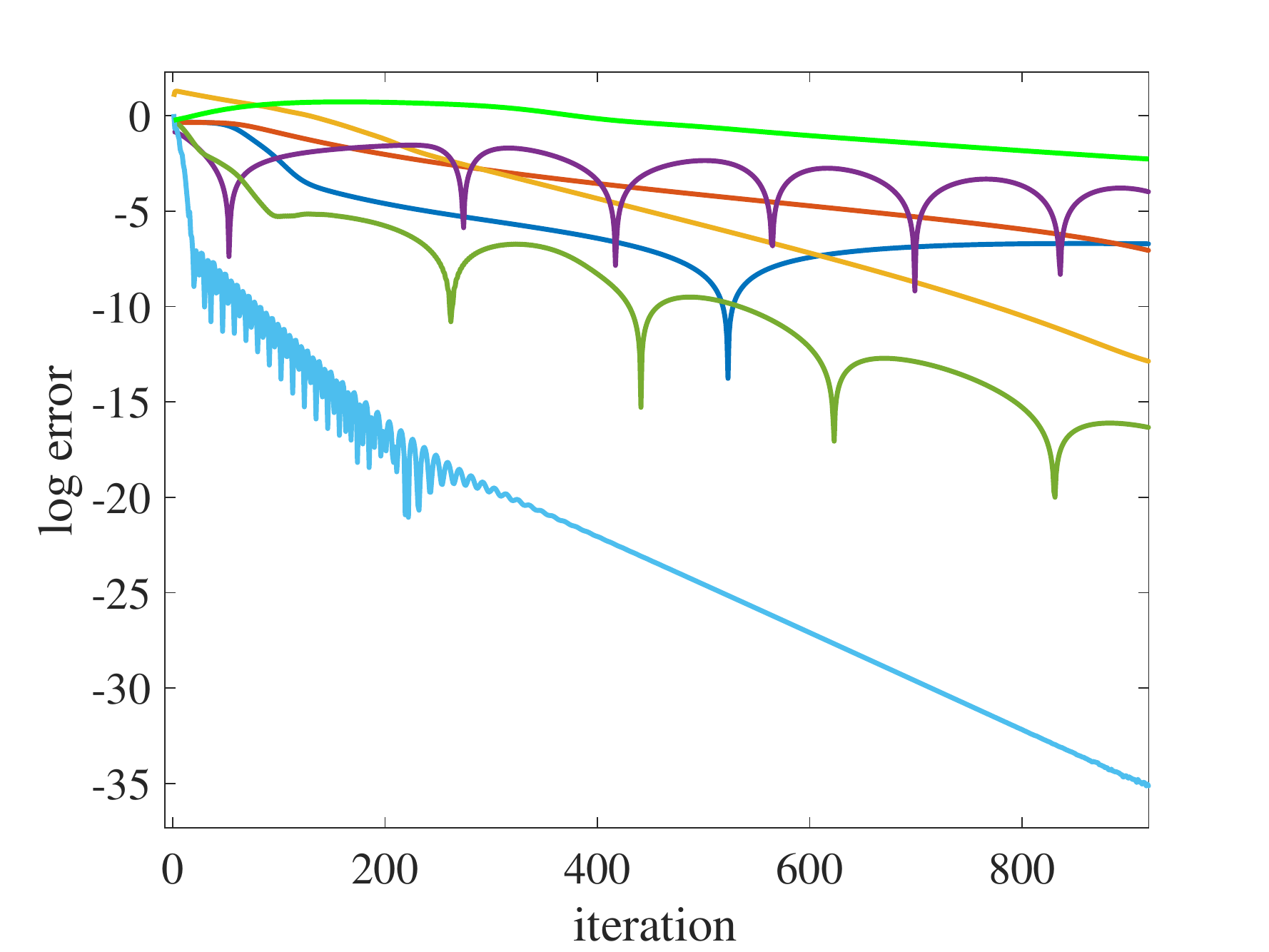}
\includegraphics[width=0.325\textwidth,trim={5 0 10 0},clip]{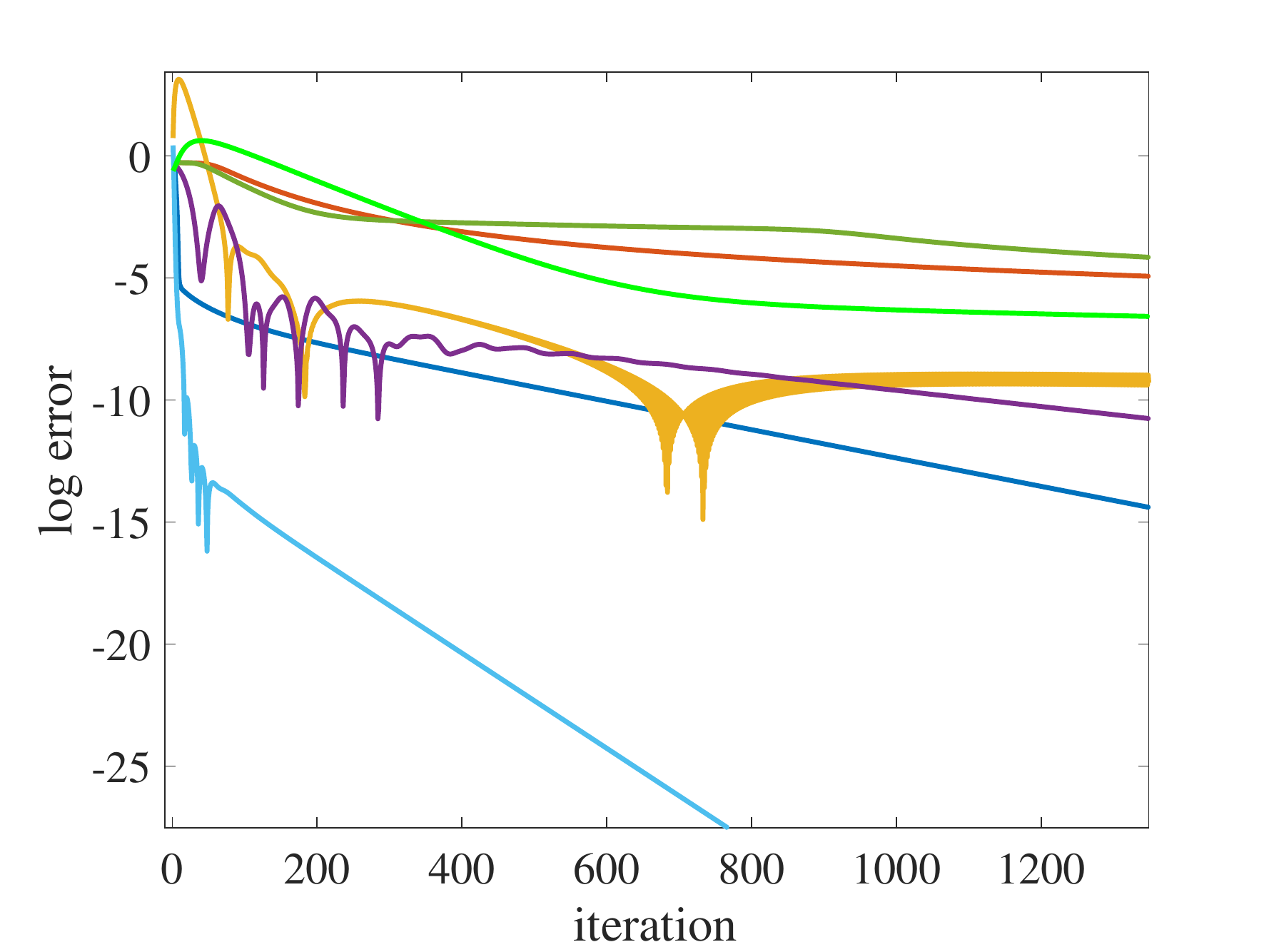}
\caption{\label{fig:nonconvexx}
Comparison between three different runs of ADMM with other algorithms when solving \eqref{eq:application_problem}.
\emph{Left:} Erd\" os-Renyi graph with edge probability $p = 3 \log(n) / n$.
\emph{Middle:} Ring with $k$-hop links with $k= \log n$.
\emph{Right:} ring graph. We plot
$\log|f(\vz_t) - f(\vz_\infty)|$ versus the iteration
$k$, where $\vx_\infty$ is the output
of the algorithm after a very large number of iterations.
}
\end{figure*}

Now  we compare the above algorithms when
solving a practical problem on sensor localization:
%
\begin{equation} \label{eq:application_problem}
    \min_{\vz_i \in \mathbb{R}^k \forall i\in\V }  \sum_{(i,j) \in \E}
    \big| \, \|\vz_i - \vz_j \|^p - d^p_{i,j}\, \big|^q + \frac{\delta}{2}\frac{|\V|}{|\E|} \sum_{i \in \V} \|\vz_i - {\bm t}\|^2,
\end{equation}
where $i\in \V$, $\{d_{i,j}\}_{i,j\in \V}$ is a set of distances between nodes, ${\bm t}$ is a constant vector,  
and $p, q > 0$.
As for the canonical problem, one can write the above problem in the same form as \eqref{eq:pair_interaction_problem}, or \eqref{eq:consensus}. Furthermore, for $d = 0$, $\delta = 0$ and $p+q= 2$ this reduces exactly to the canonical problem except that each variable node carries a vector $\vz_i$ instead of a single number $z_i$. However, in general this problem is nonconvex, and may also
be nonsmooth, e.g., when $p=q=1$.

The problem \eqref{eq:application_problem} has several important practical applications.
For instance,
given a network of $|\V|$ sensors laid out in space, where the $i$th and $j$th sensors are capable of jointly
estimating their distance $d_{i,j}$, this problem  seeks for an
accurate position for each sensor. More abstractly, given a set of distances $d$ between $|\V|$ objects, problem \eqref{eq:application_problem} finds an embedding of these objects into Euclidean space such that their
estimated distances in $\mathbb{R}^n$ are close to their true distances $d$.

When $\delta = 0$, problem \eqref{eq:application_problem}
is invariant under rotations and translations, and has an infinite number of solutions.
If $\delta \neq 0$, but ${\bm t = 0}$, it is invariant under rotations (around the origin).

All of the proximal algorithms discussed in this paper can be efficiently implemented to solve \eqref{eq:application_problem} for $p,q \in \{1,2\}$ and any $\delta$. In particular, if we assign one agent for
each term associated with each edge $(i,j) \in \E$ in the objective \eqref{eq:application_problem}, the resulting proximal maps can be computed in closed form. %
After a few changes of variables, these proximal maps are obtained by solving the one-dimensional problem
\begin{equation}\label{eq:simple_PO_for_sensor}
\min_{x \in \mathbb{R}} \big|\, |x|^p - d \, \big|^q + \frac{\rho}{2}(x - n)^2
\end{equation}
for arbitrary $d$ and $n$. For $p,q\in \{1, 2\}$, problem \eqref{eq:simple_PO_for_sensor} can be solved by finding zeros of cubic polynomials.
On the other hand, we only implemented gradient or conjugate gradient based methods for $p = q = 2$, such that each term in \eqref{eq:application_problem} is differentiable so that
all updates have closed form expressions. In particular, our conjugate gradient computations, e.g., for MSDA2, amounts to finding
the roots of cubic polynomials.
To help weaker algorithms,  we choose $\delta = 1$, ${\bm t} = {\bm 1}$, and $d_{i,j}= 1$, for all $i,j \in \V$.
%
The term multiplying $\delta$ controls the curvature of the objective
and reduces the number of local minima.


The local minimum to which different algorithms converge is strongly dependent on the initialization.
Since we are more interested in comparing convergence rates rather than the quality of local minima,
in  Fig.~\ref{fig:nonconvexx} we plot $\log( | f(\vz_t) - f(\vz_\infty) | )$ versus the iteration number $t$. All these methods are tuned with a grid search on their  parameter space. Here $\vz_\infty$ is the solution provided
by the algorithm with a very large number of iterations.
Note that the asymptotic rate of ADMM is considerably
faster than the alternatives.

\begin{figure}
\centering
\includegraphics[scale=0.35,trim={30 180 60 200},clip]{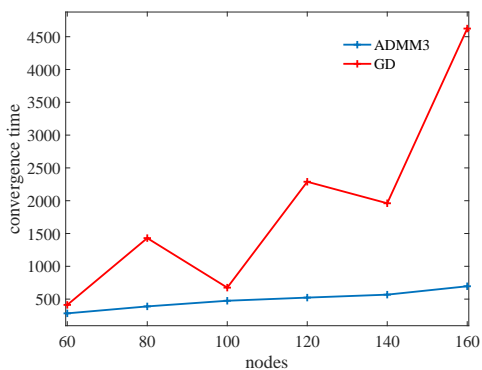}
\caption{
Ring graph. GD's convergence 
scales super linearly  with $n$ while  
ADMM scales approximately linearly.}
\label{fig:time}
\end{figure}

Finally, to check if the square root improvement of ADMM over  GD, as predicted by
Theorem~\ref{thm:proof_conj}, can be seen in practice for a problem different than the
canonical problem \eqref{eq:gen_quad}, we consider  \eqref{eq:application_problem} with a ring graph where we vary the number
of nodes and
measure the convergence time to achieve an objective value that is $\epsilon$-accurate (we choose $\epsilon \sim 10^{-20}$ for ADMM and  $\epsilon \sim 10^{-10}$ for GD).
The results are shown in Fig.~\ref{fig:time}. In this experiment, the parameters of both algorithms were  tuned with Bayesian optimization. Fig.~\ref{fig:time}  suggests that
GD scales with $n^2$ while ADMM scales with $n$, despite the nonconvexity of problem.

\section{Conclusion}

We described and compared  recent algorithms and theoretical results on distributed optimization. More specifically, we surveyed how the convergence rate of different algorithms, and  specifically of ADMM, depends
on the topology of the underlying network that constrains how different agents
communicate when solving a large consensus problem.
Since an important component of general purpose distributed  solvers consists of a distributed averaging subroutine, we also focused on
a particular
distributed averaging problem, namely  \eqref{eq:gen_quad}. Regarding ADMM, we showed an explicit, and optimal, rate of convergence   analysis for this problem.
We related the optimal convergence rate
of ADMM with the second largest eigenvalue of the transition matrix
of the communicating  network, and also provided explicit formulae for optimal parameter tuning in terms of this
eigenvalue.
We also showed that ADMM can be seen as a \emph{lifting} of GD, in close analogy
with lifted Markov chains theory. We showed that ADMM attains a square root
speedup over GD that is reminiscent of the maximum possible mixing achieved via lifted Markov chains, but which however is only possible for some types of graph with not so small conductance. On the other hand, in the case of ADMM, such a relation
holds for any graph. These results
provide interesting connections
between distributed optimization and other fields
of mathematics and may be of independent interest.

We also verified numerically that our theoretical results match the practice. We numerically
compared ADMM with several state-of-the-art methods regarding the performance of the distributed averaging subroutines (via \eqref{eq:gen_quad}) and the performance of these methods on a related nonconvex problem.
ADMM proved to be in general faster than competing methods.

\bibliographystyle{unsrt}
\bibliography{biblio.bib}

\end{document}